\documentclass{siamart251216}

\usepackage{amsfonts}
\usepackage{amssymb}
\usepackage{mathtools}
\usepackage{commath}
\usepackage{algorithm}
\usepackage{algorithmic}
\usepackage{graphicx}
\usepackage[caption=false]{subfig}
\usepackage{textcomp}
\usepackage{nicefrac}
\usepackage{bm}
\usepackage{bbm}
\usepackage{makecell}
\usepackage{xcolor}
\usepackage{tikz}
\usetikzlibrary{arrows.meta, positioning}

\ifpdf
  \DeclareGraphicsExtensions{.eps,.pdf,.png,.jpg}
\else
  \DeclareGraphicsExtensions{.eps}
\fi

\newsiamthm{assumption}{Assumption}
\newsiamremark{remark}{Remark}

\headers{Noisy Gradient for Resource Allocation}{L. Qin and Y. Pu}

\title{On the Convergence of a Noisy Gradient Method for Non-convex Distributed Resource Allocation: Saddle Point Escape\thanks{This work was supported by a Melbourne Research Scholarship and the Australian Research Council (DE220101527).}}

\author{
Lei Qin
\thanks{Department of Electrical and Electronic Engineering, University of Melbourne,
Parkville VIC 3010, Australia
(\email{leqin@student.unimelb.edu.au}, \email{ye.pu@unimelb.edu.au}).}
\and
Ye Pu\footnotemark[2]
}

\ifpdf
\hypersetup{
  pdftitle={On the Convergence of a Noisy Gradient Method for Non-convex Distributed Resource Allocation: Saddle Point Escape},
  pdfauthor={Lei Qin and Ye Pu}
}
\fi

\begin{document}
\maketitle

\begin{abstract}
This paper considers a class of distributed resource allocation problems where each agent privately holds a smooth, potentially non-convex local objective, subject to a globally coupled equality constraint. Built upon the existing method, Laplacian-weighted Gradient Descent, we propose to add random perturbations to the gradient iteration to enable efficient escape from saddle points and achieve second-order convergence guarantees. We show that, with a sufficiently small fixed step size, the iterates of all agents reach an approximate second-order stationary point with high probability. Numerical experiments confirm the effectiveness of the proposed approach, demonstrating improved performance over standard weighted gradient descent in non-convex scenarios.
\end{abstract}



\section{Introduction}
\label{sec: Introduction}
Distributed resource allocation aims to optimize local objectives over a network while satisfying a global resource constraint, with applications including economic dispatch \cite{li2021distributed,yi2016initialization,wang2018distributed}, communication-resource allocation \cite{halabian2019distributed}, and multi-agent task allocation \cite{sayyaadi2010distributed}; see \cite{doostmohammadian2025survey} for a recent survey.

In particular, we consider a resource allocation problem over a network of $m$ agents subject to the global resource constraint
\begin{gather}
    \label{eq: Resource allocation problem}
    \begin{aligned}
        \min_{\boldsymbol{\theta} \in (\mathbb{R}^{n})^{m}} \quad & F(\boldsymbol{\theta}) \triangleq \sum_{i=1}^{m} f_{i}(\boldsymbol{\theta}_{i})\\
        \text{subject to} \quad & \sum_{i=1}^{m} \boldsymbol{\theta}_{i} = \mathbf{r},
    \end{aligned}
\end{gather}
where $\boldsymbol{\theta}=[(\boldsymbol{\theta}_{1})^{\top},\dots,(\boldsymbol{\theta}_{m})^{\top}]^{\top}\in(\mathbb{R}^{n})^{m}$ is the decision vector and $\mathbf{r}\in\mathbb{R}^{n}$ is a given resource demand vector. Each local objective $f_i:\mathbb{R}^{n}\rightarrow\mathbb{R}$ is smooth, possibly non-convex, and privately known only to agent $i$. The network is modeled as an undirected and connected graph $\mathcal{G}(\mathcal{V},\mathcal{E})$ with node set $\mathcal{V}:=\{1,\ldots,m\}$ and edge set $\mathcal{E}\subseteq\mathcal{V}\times\mathcal{V}$. Agent $i$ can communicate directly with agent $j$ if $(i,j)\in\mathcal{E}$.

Distributed resource allocation has been extensively studied in convex settings, including ADMM-based methods \cite{chang2014multi,chang2016proximal,aybat2019distributed}, dual and consensus approaches \cite{zhang2019distributed,yuan2019stochastic}, first-order and accelerated methods \cite{beck20141,necoara2013random,ghadimi2011accelerated}, and continuous-time schemes \cite{zhu2019distributed,deng2017distributed}. Among these approaches, \textbf{Laplacian-weighted Gradient Descent (LGD)} \cite{xiao2006optimal,ho1980class} is particularly attractive because of its simple structure and guaranteed feasibility at every iteration. Its fixed-step update is
\begin{gather}
    \label{eq: LGD}
    \boldsymbol{\theta}_{i}^{k+1}=\boldsymbol{\theta}_{i}^{k}-\alpha\sum_{j=1}^{m}\ell_{ij}\nabla f_{j}(\boldsymbol{\theta}_{j}^{k}),
\end{gather}
where $\alpha>0$ is the step size and $\ell_{ij}$ is the $(i,j)$-th entry of the Laplacian matrix $\mathbf{L}$. Related developments address delays and communication constraints \cite{doostmohammadian20221st,doostmohammadian2022distributed}, acceleration and momentum \cite{ghadimi2011accelerated,ochoa2019hybrid,doostmohammadian2024accelerated}, and feasibility correction \cite{li2022implicit}.

For non-convex resource allocation, distributed approaches include generalized Lagrangian formulations and gradient-tracking methods \cite{li2020projection,zhang2020distributed}, as well as privacy-preserving formulations \cite{Beaude2020}. Related distributed non-convex problems with separable objectives and coupled affine constraints have been studied using augmented Lagrangian methods \cite{Houska2016}, while limitations of decentralized non-convex optimization are analyzed in \cite{Hong2022}. A momentum-based resource-allocation method under additional second-order assumptions is proposed in \cite{xia2025distributed}. Most existing methods in this setting provide first-order guarantees and therefore do not, in general, exclude strict saddle points. Moreover, gradient methods can require exponentially many iterations to escape saddle points in the worst case \cite{du2017gradient}.

Our goal is to obtain second-order guarantees using only first-order information. Random perturbations are a standard mechanism for escaping strict saddle points in centralized non-convex optimization \cite{ge2015escaping,jin2017escape,jin2021nonconvex,Davis2022}. Second-order properties of distributed gradient and gradient-tracking methods have also been studied in \cite{Daneshmand2020}, while perturbation-based distributed methods are considered in \cite{vlaski2021distributed,wang2023decentralized,qin2023second}. These works, however, primarily concern unconstrained distributed or consensus optimization. In contrast, we consider distributed resource allocation with a global coupling constraint. We structure the perturbations through the square-root Laplacian so that they lie in the disagreement subspace, preserving feasibility at every iteration while enabling escape from strict saddle points.

The main contributions are summarized as follows:
\begin{itemize}
    \item We develop an equivalent representation of \textbf{LGD} for distributed resource allocation (Proposition~\ref{pro: Same sequence}), which enables us to characterize its first-order stationarity properties for non-convex objectives (Proposition~\ref{pro: First order property}).
    \item We propose \textbf{Noisy LGD (NLGD)}, which augments \textbf{LGD} with structured random perturbations to obtain second-order guarantees.
    \item We prove that \textbf{NLGD} reaches an approximate second-order stationary point with high probability (Theorem~\ref{the: Approximate second order}), while preserving the global resource constraint at every iteration and hence ensuring \emph{anytime feasibility}.
\end{itemize}

The remainder of the paper presents the assumptions and supporting results in Section~\ref{sec: Assumptions and Supporting Results}, the proposed algorithm and main results in Section~\ref{sec: Main Results}, the proofs in Section~\ref{sec: Proofs}, and numerical examples in Section~\ref{sec: Numerical Examples}.

\subsection{Notation}
\label{sec: Notation}
Let $\mathbf{I}_n$ denote the $n\times n$ identity matrix, $\bm{1}_n$ the $n$-vector of ones, and $a_{ij}$ the $(i,j)$-th entry of a matrix $\mathbf{A}$. For a symmetric matrix $\mathbf{B}$, $\lambda_{\min}(\mathbf{B})$, $\lambda_{\max}(\mathbf{B})$, and $\|\mathbf{B}\|$ denote its minimum eigenvalue, maximum eigenvalue, and spectral norm, respectively. For a symmetric positive semidefinite matrix $\mathbf{C}$, $\lambda_{\min}^{+}(\mathbf{C})$ denotes its smallest nonzero eigenvalue. The Kronecker product is denoted by $\otimes$. We write $\mathbf{x}\sim\mathcal{N}(\bm{\mu},\bm{\Sigma})$ for an $n$-dimensional Gaussian random vector with mean $\bm{\mu}\in\mathbb{R}^{n}$ and covariance matrix $\bm{\Sigma}\in\mathbb{R}^{n\times n}$. The ceiling function is denoted by $\lceil\cdot\rceil$. Unless explicitly stated otherwise, iteration indices are nonnegative integers.

\section{Assumptions and Supporting Results}
\label{sec: Assumptions and Supporting Results}
\subsection{Assumptions}
\begin{assumption}[Lipschitz continuity]
\label{ass: Lipschitz}
    Each $f_{i}$ in \eqref{eq: Resource allocation problem} is both $L_{f_{i}}^{g}$-gradient Lipschitz and $L_{f_{i}}^{H}$-Hessian Lipschitz, i.e., for all $\boldsymbol{\theta}_{i},\boldsymbol{\omega}_{i}\in\mathbb{R}^{n}$ and each $i \in \mathcal{V}$,
    \begin{align*}
        \|\nabla f_{i}(\boldsymbol{\theta}_{i}) - \nabla f_{i}(\boldsymbol{\omega}_{i})\|
        &\le L_{f_{i}}^{g} \|\boldsymbol{\theta}_{i}-\boldsymbol{\omega}_{i}\|,\\
        \|\nabla^{2} f_{i}(\boldsymbol{\theta}_{i}) - \nabla^{2} f_{i}(\boldsymbol{\omega}_{i})\|
        &\le L_{f_{i}}^{H} \|\boldsymbol{\theta}_{i}-\boldsymbol{\omega}_{i}\|.
    \end{align*}
\end{assumption}

\begin{assumption}[Coercivity]
\label{ass: Coercivity}
    Each $f_{i}$ in \eqref{eq: Resource allocation problem} is coercive (i.e., its sublevel sets are compact by continuity).
\end{assumption}

\begin{remark}
    If Assumption \ref{ass: Lipschitz} holds, then $F$ defined in \eqref{eq: Resource allocation problem} has $L_{F}^{g}$-Lipschitz continuous gradient and $L_{F}^{H}$-Lipschitz continuous Hessian with
    \begin{gather}
    \label{eq: Lipschitz continuous for F}
        L_{F}^{g} = \max_{i} \{L_{f_{i}}^{g}\},~L_{F}^{H} = \max_{i} \{L_{f_{i}}^{H}\}.
    \end{gather}
    If Assumption \ref{ass: Coercivity} holds, then $F$ defined in \eqref{eq: Resource allocation problem} is also coercive.
\end{remark}

\begin{assumption}[Connected network]
\label{ass: Network}
    The undirected network graph $\mathcal{G}(\mathcal{V},\mathcal{E})$ is connected.
\end{assumption}

\subsection{Laplacian-weighted Gradient Descent}
In this section, we review the standard weighted gradient method in \eqref{eq: LGD} for distributed resource allocation and study its convergence properties. First, recall the objective function $F(\boldsymbol{\theta})$ defined in \eqref{eq: Resource allocation problem}. Note that,
\begin{align*}
    \nabla F(\boldsymbol{\theta})
    &= [\nabla f_{1}(\boldsymbol{\theta}_{1})^{\top},\cdots, \nabla f_{m}(\boldsymbol{\theta}_{m})^{\top}]^{\top},\\
    \nabla^{2} F(\boldsymbol{\theta})
    &= \bigoplus_{i=1}^{m} \nabla^{2} f_{i}(\boldsymbol{\theta}_{i}),
\end{align*}
where $\bigoplus$ denotes the block diagonal concatenation of matrices. In particular, the Hessian of $F$ is block diagonal.

The fixed step-size \textbf{LGD} in \eqref{eq: LGD} can be formulated in an aggregate form as
\begin{gather}
\label{eq: LGD aggregate form}
    \boldsymbol{\theta}^{k+1} = \boldsymbol{\theta}^{k} - \alpha \hat{\mathbf{L}} \cdot \nabla  F(\boldsymbol{\theta}^{k})
\end{gather}
with $\hat{\mathbf{L}} = \mathbf{L}\otimes \mathbf{I}_{n}$.

\begin{remark}
    The matrix $\mathbf{L}$ can be replaced by any symmetric weighting matrix with zero row sums. Without loss of generality, we use the Laplacian matrix. Since $\hat{\mathbf{L}} = \mathbf{L}\otimes \mathbf{I}_{n}$ and $\mathbf{L}$ is a symmetric positive semi-definite matrix, there exists a unique symmetric positive semi-definite matrix $\sqrt{\hat{\mathbf{L}}} \in (\mathbb{R}^{n})^{m}$ such that $\sqrt{\hat{\mathbf{L}}} \cdot \sqrt{\hat{\mathbf{L}}} = \hat{\mathbf{L}}$. 
\end{remark}

The following result shows that applying \textbf{LGD} to Problem \eqref{eq: Resource allocation problem} is equivalent to performing gradient descent on the auxiliary function $\Psi_{\boldsymbol{\theta}^{0}}$ with its proof provided in Section \ref{sec: Proofs}.
\begin{proposition}
\label{pro: Same sequence}
    Let auxiliary function
    \begin{gather}
        \begin{gathered}
            \label{eq: Auxiliary function}
            \Psi_{\boldsymbol{\theta}^{0}}(\mathbf{x}) \triangleq F(\boldsymbol{\theta}^{0} + \sqrt{\hat{\mathbf{L}}} \mathbf{x}).
        \end{gathered}    
    \end{gather}
    For the distributed resource allocation problem in \eqref{eq: Resource allocation problem}, given fixed step-size $\alpha > 0$ and initial point $\boldsymbol{\theta}^{0} \in (\mathbb{R}^{n})^{m}$ satisfying $\bm{1}_{m}^{\top} \otimes \mathbf{I}_{n} \cdot \boldsymbol{\theta}^{0} = \bm{r}$, the sequence $\{\boldsymbol{\theta}^{k}\}$ generated by \eqref{eq: LGD aggregate form} is equivalent to the sequence generated by applying gradient descent to $\Psi_{\boldsymbol{\theta}^{0}}$ from the same initial point $\boldsymbol{\theta}^{0} \in (\mathbb{R}^{n})^{m}$ and $\mathbf{x}^{0} = \bm{0} \in (\mathbb{R}^{n})^{m}$ with the same fixed step-size $\alpha > 0$, as per
    \begin{gather}
    \label{eq: LGD with auxiliary function}
        \begin{gathered}
            \mathbf{x}^{k+1} = \mathbf{x}^{k} - \alpha \nabla \Psi_{\boldsymbol{\theta}^{0}}(\mathbf{x}^{k}),\\
            \boldsymbol{\theta}^{k+1} = \boldsymbol{\theta}^{0} + \sqrt{\hat{\mathbf{L}}} \mathbf{x}^{k+1}.
        \end{gathered}
    \end{gather}
\end{proposition}

\begin{remark}
    If Assumption \ref{ass: Lipschitz} holds, then $\Psi_{\boldsymbol{\theta}^{0}}$ in \eqref{eq: Auxiliary function} has $L_{\Psi_{\boldsymbol{\theta}^{0}}}^{g}$-Lipschitz continuous gradient and $L_{\Psi_{\boldsymbol{\theta}^{0}}}^{H}$-Lipschitz continuous Hessian with
    \begin{gather}
    \label{eq: Psi Lipschitz}
        L_{\Psi_{\boldsymbol{\theta}^{0}}}^{g} = \norm{\sqrt{\mathbf{L}}}^{2} \cdot L_{F}^{g},~L_{\Psi_{\boldsymbol{\theta}^{0}}}^{H} =  \norm{\sqrt{\mathbf{L}}}^{3} \cdot L_{F}^{H}.
    \end{gather}
\end{remark}

\begin{definition}[adapted from Lemma 1 \cite{nedic2018improved}]
\label{def: first-order stationary point}
    For the distributed resource allocation problem in \eqref{eq: Resource allocation problem}, a point $\boldsymbol{\theta} \in (\mathbb{R}^{n})^{m}$ is said to be a first-order stationary point if it satisfies the following:
    \begin{enumerate}
    \renewcommand{\labelenumi}{\roman{enumi})}
        \item $\sqrt{\hat{\mathbf{L}}} \cdot \nabla F(\boldsymbol{\theta}) = 0$;
        \item $\bm{1}_{m}^{\top} \otimes \mathbf{I}_{n} \cdot \boldsymbol{\theta} = \bm{r}$,
    \end{enumerate}
    where $\sqrt{\hat{\mathbf{L}}} = \sqrt{\mathbf{L}}\otimes \mathbf{I}_{n}$ with $\sqrt{\mathbf{L}} \cdot \sqrt{\mathbf{L}} = \mathbf{L}$.
\end{definition}

\begin{remark}
    If condition i) in Definition \ref{def: first-order stationary point} holds, then $\nabla f_{i}(\boldsymbol{\theta}_{i})$ is in consensus, i.e., $\nabla f_{i}(\boldsymbol{\theta}_{i}) = \nabla f_{j}(\boldsymbol{\theta}_{j})$ for all $i$, $j \in \mathcal{V}$. Furthermore, if condition ii) also holds, then $\boldsymbol{\theta}$ satisfies the first-order stationarity conditions of Problem \eqref{eq: Resource allocation problem}.
\end{remark}
 
The following result, with its proof provided in Section \ref{sec: Proofs}, shows the first order stationarity guarantees of \textbf{LGD} update \eqref{eq: LGD aggregate form}, and its proof is based on $\Psi_{\boldsymbol{\theta}^{0}}$ defined in \eqref{eq: Auxiliary function}.

\begin{proposition}
\label{pro: First order property}
    Let Assumptions \ref{ass: Lipschitz}, \ref{ass: Coercivity} and \ref{ass: Network} hold. Given initial point $\boldsymbol{\theta}^{0} \in (\mathbb{R}^{n})^{m}$ satisfying $\bm{1}_{m}^{\top} \otimes \mathbf{I}_{n} \cdot \boldsymbol{\theta}^{0} = \bm{r}$, for any fixed step-size
    \begin{gather*}
        0 < \alpha \le \frac{1}{\|\sqrt{\mathbf{L}}\|^{2} \cdot L_{F}^{g}},
    \end{gather*}
    the sequence $\{\boldsymbol{\theta}^{k}\}$ generated by \eqref{eq: LGD aggregate form} satisfies that for all $k \ge 0$, $\bm{1}_{m}^{\top} \otimes \mathbf{I}_{n} \cdot \boldsymbol{\theta}^{k}= \bm{r}$ and
    \begin{gather*}
        \lim_{k\rightarrow\infty} \norm{\sqrt{\hat{\mathbf{L}}} \cdot \nabla F(\boldsymbol{\theta}^{k})} = 0.
    \end{gather*}
\end{proposition}

Next, we introduce the definition of an approximately consensual second-order stationary point.
\begin{definition}
\label{def: Approximately second order stationary point}
    For the distributed resource allocation problem in \eqref{eq: Resource allocation problem}, a point $\boldsymbol{\theta} \in (\mathbb{R}^{n})^{m}$ is said to be an $(\epsilon,\gamma)$-second-order stationary point if it satisfies the following:
    \begin{enumerate}
    \renewcommand{\labelenumi}{\roman{enumi})}
        \item $\|\sqrt{\hat{\mathbf{L}}} \cdot \nabla F(\boldsymbol{\theta})\| \le \epsilon$;
        \item $\bm{1}_{m}^{\top} \otimes \mathbf{I}_{n} \cdot \boldsymbol{\theta} = \bm{r}$;
        \item $\mathbf{d}^{\top} \nabla^{2} F(\boldsymbol{\theta}) \mathbf{d} \ge -\gamma \norm{\mathbf{d}}^{2}$ for all $\mathbf{d} \in \mathcal{T}$,
    \end{enumerate} 
    where $\sqrt{\hat{\mathbf{L}}} = \sqrt{\mathbf{L}}\otimes \mathbf{I}_{n}$ with $\sqrt{\mathbf{L}} \cdot \sqrt{\mathbf{L}} = \mathbf{L}$ and tangent space $\mathcal{T} = \{\mathbf{d} \in (\mathbb{R}^{n})^{m}~:~\bm{1}_{m}^{\top} \otimes \mathbf{I}_{n} \cdot \mathbf{d} = \bm{0}\}$.
\end{definition}

Condition i) and ii) in Definition \ref{def: Approximately second order stationary point} mean $\boldsymbol{\theta}$ is an approximate first-order stationary point by Definition \ref{def: first-order stationary point}. Further, if condition iii) holds, $\nabla^{2} F(\boldsymbol{\theta})$ is not excessively negative on the orthogonal complement of $\text{span}\{\bm{1}_{m}^{\top} \otimes \mathbf{I}_{n}\}$, i.e., the feasible directions. We generically refer to such points as approximately second-order stationary points. These approximate second-order stationary points include local minimizers and exclude saddle points with significant negative curvature. The following result establishes  a connection between the approximate second-order stationary points of the auxiliary function $\Psi_{\boldsymbol{\theta}^{0}}$ and the approximate second-order stationary points of the original objective $F$ in Problem \eqref{eq: Resource allocation problem}. The corresponding proof is provided in Section \ref{sec: Proofs}.

\begin{proposition}
\label{pro: Second order relations}
    Let Assumption \ref{ass: Network} hold. For the distributed resource allocation problem in \eqref{eq: Resource allocation problem}, given an initial point $\boldsymbol{\theta}^{0} \in (\mathbb{R}^{n})^{m}$ satisfying $\bm{1}_{m}^{\top} \otimes \mathbf{I}_{n} \cdot \boldsymbol{\theta}^{0} = \bm{r}$, if the following holds at $\mathbf{x} \in (\mathbb{R}^{n})^{m}$:
    \begin{enumerate}
    \renewcommand{\labelenumi}{\roman{enumi})}
        \item $\|\nabla \Psi_{\boldsymbol{\theta}^{0}}(\mathbf{x})\| \le \epsilon$;
        \item $\lambda_{\min} (\nabla^{2} \Psi_{\boldsymbol{\theta}^{0}}(\mathbf{x})) \ge -\gamma$,
    \end{enumerate}
    then, $\boldsymbol{\theta} = \boldsymbol{\theta}^{0} + \sqrt{\hat{\mathbf{L}}} \mathbf{x}$ is an $(\epsilon,\gamma/\lambda_{\min}^{+}(\mathbf{L}))$-second-order stationary point of Problem \eqref{eq: Resource allocation problem}, where $\Psi_{\boldsymbol{\theta}^{0}}$ is defined in \eqref{eq: Auxiliary function} and $\sqrt{\hat{\mathbf{L}}} = \sqrt{\mathbf{L}}\otimes \mathbf{I}_{n}$ with $\sqrt{\mathbf{L}} \cdot \sqrt{\mathbf{L}} = \mathbf{L}$.
\end{proposition}

\section{Main Results}
\label{sec: Main Results}
To obtain second-order guarantees using only first-order information, we propose \textbf{Noisy LGD (NLGD)}. Given a feasible reference point $\boldsymbol{\theta}^{0}$, each agent $i\in\mathcal{V}$ updates
\begin{gather}
    \label{eq: NLGD}
    \boldsymbol{\theta}_{i}^{k+1}=\boldsymbol{\theta}_{i}^{k}-\alpha\sum_{j=1}^{m}\left(\ell_{ij}\nabla f_j(\boldsymbol{\theta}_{j}^{k})+\tilde{\ell}_{ij}\mathbf{n}_{j}^{k}\right),
\end{gather}
where $\ell_{ij}$ and $\tilde{\ell}_{ij}$ are the $(i,j)$-th entries of $\mathbf{L}$ and its symmetric square root $\sqrt{\mathbf{L}}$, respectively. The perturbations $\mathbf{n}_{j}^{k}$ are generated from a common pseudo-random generator with shared seed $s$, so that every agent can reproduce the required perturbation vectors locally without additional communication. The matrices $\mathbf{L}$ and $\sqrt{\mathbf{L}}$, the step size $\alpha$, and the noise variance $\sigma^{2}$ are selected offline.

\begin{algorithm}[h]
    \caption{Noisy Laplacian Gradient Descent (NLGD)}
    \label{alg:NLGD}
    \begin{algorithmic}[1]
    \REQUIRE Feasible $\boldsymbol{\theta}_{i}^{0}$, $\alpha$, $\{\ell_{ij},\tilde{\ell}_{ij}\}_{j=1}^{m}$, $\sigma^{2}$, and shared seed $s$
    \FOR{$k = 0,1,2,\dots$}
        \FORALL{agents $i \in \mathcal{V}$}
            \STATE Exchange $\nabla f_i(\boldsymbol{\theta}_{i}^{k})$ with neighboring agents $j$ satisfying $\ell_{ij}\neq 0$
            \STATE Generate perturbations $\mathbf{n}_j^k = \mathrm{PRNG}(s,k,j)$ for all $j \in \mathcal{V}$ using a common pseudo-random generator with shared seed $s$
            \STATE Compute $\boldsymbol{\theta}_{i}^{k+1}$ using \eqref{eq: NLGD}
        \ENDFOR
    \ENDFOR
    \end{algorithmic}
\end{algorithm}

\begin{remark}
    The square root $\sqrt{\mathbf{L}}$ can be computed offline when the topology is known. Alternatively, $\sqrt{\mathbf{L}}\mathbf{n}^{k}$ can be approximated using distributed polynomial filtering, trading additional communication rounds for approximation error.
\end{remark}

In aggregate form, \eqref{eq: NLGD} becomes
\begin{gather}
\label{eq: NLGD aggregate form}
    \begin{gathered}
       \boldsymbol{\theta}^{k+1} = \boldsymbol{\theta}^{k} - \alpha (\hat{\mathbf{L}} \cdot\nabla  F(\boldsymbol{\theta}^{k}) + \sqrt{\hat{\mathbf{L}}}\mathbf{n}^{k}),    
    \end{gathered}
\end{gather}
where $\mathbf{n}^{k}=[(\mathbf{n}_{1}^{k})^{\top},\ldots,(\mathbf{n}_{m}^{k})^{\top}]^{\top}\in(\mathbb{R}^{n})^{m}$.

The following result extends Proposition~\ref{pro: Same sequence} to \textbf{NLGD}; the proof is identical and omitted.
\begin{proposition}
\label{pro: NGD interpretation}
Let $\Psi_{\boldsymbol{\theta}^{0}}$ be defined in \eqref{eq: Auxiliary function}. For any $\alpha>0$, feasible $\boldsymbol{\theta}^{0}$ satisfying $\bm{1}_{m}^{\top}\otimes\mathbf{I}_{n}\cdot\boldsymbol{\theta}^{0}=\mathbf{r}$, and perturbation sequence $\{\mathbf{n}^{k}\}$, the sequence $\{\boldsymbol{\theta}^{k}\}$ generated by \eqref{eq: NLGD aggregate form} is equivalently generated from $\mathbf{x}^{0}=\bm{0}$ by
\begin{subequations}
\label{eq: NLGD with auxiliary function}
\begin{align}
\mathbf{x}^{k+1}&=\mathbf{x}^{k}-\alpha\left(\nabla\Psi_{\boldsymbol{\theta}^{0}}(\mathbf{x}^{k})+\mathbf{n}^{k}\right), \label{eq: NLGD with auxiliary function 1}\\
\boldsymbol{\theta}^{k+1}&=\boldsymbol{\theta}^{0}+\sqrt{\hat{\mathbf{L}}}\mathbf{x}^{k+1}. \label{eq: NLGD with auxiliary function 2}
\end{align}
\end{subequations}
\end{proposition}

\begin{assumption}[Random perturbation]
\label{ass: Random perturbation}
The perturbations $\{\mathbf{n}_{i}^{k}\}_{i,k}$ are independent and satisfy $\mathbf{n}_{i}^{k}\sim\mathcal{N}(\bm{0},\sigma^{2}\mathbf{I}_{n})$ for all $k\geq0$ and $i\in\mathcal{V}$. Hence, $\mathbf{n}^{k}=[(\mathbf{n}_{1}^{k})^{\top},\ldots,(\mathbf{n}_{m}^{k})^{\top}]^{\top}\sim\mathcal{N}(\bm{0},\sigma^{2}\mathbf{I}_{mn})$ and $\mathbb{E}[\|\mathbf{n}^{k}\|^{2}]=mn\sigma^{2}$.
\end{assumption}

We first establish the second-order guarantee for noisy gradient descent applied to the auxiliary function $\Psi_{\boldsymbol{\theta}^{0}}$.
\begin{proposition}
\label{pro: Approximate second order}
    Let Assumptions \ref{ass: Lipschitz}, \ref{ass: Coercivity}, \ref{ass: Network} and \ref{ass: Random perturbation} hold. Further, let $f_{i}^{\star}$ denote the global minimum of function $f_{i}$ for $i\in\mathcal{V}$. Then, given parameter $\epsilon_{g} > 0$, $\epsilon_{H} = \sqrt{\epsilon_{g}L_{\Psi_{\boldsymbol{\theta}^{0}}}^{H}}$, and confidence parameter $0 < p < 1$ with $L_{\Psi_{\boldsymbol{\theta}^{0}}}^{g},L_{\Psi_{\boldsymbol{\theta}^{0}}}^{H}$ as per \eqref{eq: Psi Lipschitz}, there exists
    \begin{gather}
    \label{eq: alpha bar}
        \bar{\alpha} \le \min\{\frac{1}{L_{\Psi_{\boldsymbol{\theta}^{0}}}^{g}},~-\frac{2\ln{(p)}}{L_{\Psi_{\boldsymbol{\theta}^{0}}}^{g}}\}
    \end{gather}
    such that for any step-size $\alpha \le \bar{\alpha}$, with random perturbation variance
    \begin{gather}
    \label{eq: Variance}
        \sigma^{2}=\frac{L_{\Psi_{\boldsymbol{\theta}^{0}}}^{g}\alpha\epsilon_{g}^{2}}{12mn},
    \end{gather}
    and initial condition satisfying $\bm{1}_{m}^{\top} \otimes \mathbf{I}_{n} \cdot \boldsymbol{\theta}^{0} = \bm{r}$ and $\mathbf{x}^{0} = \bm{0}$, after
    \begin{gather}
    \label{eq: Iteration}
        K = \lceil \frac{\Psi_{\boldsymbol{\theta}^{0}}(\mathbf{x}^{0}) - \sum_{i=1}^{m} f_{i}^{\star}}{L_{\Psi_{\boldsymbol{\theta}^{0}}}^{g} \epsilon_{g}^{2} \alpha^{2}} \rceil
    \end{gather}
    iterations of \eqref{eq: NLGD with auxiliary function}, it follows that
    \begin{gather}
    \label{eq: Approximate second order}
        \mathbb{P} \Bigg[\exists k \in [0,K],~\|\nabla \Psi_{\boldsymbol{\theta}^{0}}(\mathbf{x}^{k})\| \le\epsilon_{g}~~\land~~\lambda_{\min} (\nabla^{2} \Psi_{\boldsymbol{\theta}^{0}}(\mathbf{x}^{k})) \ge -\epsilon_{H}
        \Bigg] \ge 1 - p.
    \end{gather}
\end{proposition}

\begin{remark}
\label{re: Discussion}
    The tolerances $\epsilon_{g}$ and $\epsilon_{H}$ quantify first- and second-order stationarity, respectively, while the iteration bound in \eqref{eq: Iteration} scales as $O(\epsilon_{g}^{-2}\alpha^{-2})$. The probability bound $1-p$ in \eqref{eq: Approximate second order} gives the confidence of reaching such a point within $K$ iterations; a higher confidence level may require a more restrictive step-size selection, as specified in the proof of Proposition~\ref{pro: Approximate second order}.
\end{remark}

The following theorem gives the main second-order guarantee for \textbf{NLGD}.

\begin{theorem}
\label{the: Approximate second order}
    Let Assumptions \ref{ass: Lipschitz}, \ref{ass: Coercivity}, \ref{ass: Network} and \ref{ass: Random perturbation} hold. Further, let $f_{i}^{\star}$ denote the global minimum of function $f_{i}$ for $i\in\mathcal{V}$. Then, given parameters $\epsilon_{g} > 0$, and $\epsilon_{H} = \sqrt{\epsilon_{g} \cdot\|\sqrt{\mathbf{L}}\|^{3} \cdot L_{F}^{H}}$, and confidence parameter $0 < p < 1$ with $L_{F}^{g},L_{F}^{H}$ as per \eqref{eq: Lipschitz continuous for F}, there exists
    \begin{gather*}
        \bar{\alpha} \le \min\{\frac{1}{\|\sqrt{\mathbf{L}}\|^{2} \cdot L_{F}^{g}},~-\frac{2\ln{(p)}}{\|\sqrt{\mathbf{L}}\|^{2} \cdot L_{F}^{g}}\}
    \end{gather*}
    such that for any step-size $\alpha \le \bar{\alpha}$, with random perturbation
    variance $\sigma^{2}$ as per \eqref{eq: Variance}, $K$ as per \eqref{eq: Iteration}, $\{\boldsymbol{\theta}^{k}\}$ given by \eqref{eq: NLGD} with an initial point $\boldsymbol{\theta}^{0} \in (\mathbb{R}^{n})^{m}$ satisfying $\bm{1}_{m}^{\top} \otimes \mathbf{I}_{n} \cdot \boldsymbol{\theta}^{0} = \bm{r}$, within $K$ iterations of \eqref{eq: NLGD}, with probability $1-p$, there exists an $(\epsilon_{g},\epsilon_{H}/\lambda_{\min}^{+}(\mathbf{L}))$-second-order stationary point of Problem \eqref{eq: Resource allocation problem}, i.e.,
    \begin{multline*}
        \mathbb{P} \Bigg[\exists k \in [0,K],~\|\sqrt{\hat{\mathbf{L}}} \cdot \nabla F(\boldsymbol{\theta}^{k})\| \le\epsilon_{g}\\
        \land~~\forall~\mathbf{d} \in \mathcal{T},~\mathbf{d}^{\top} \nabla^{2} F(\boldsymbol{\theta}^{k}) \mathbf{d} \ge - \frac{\epsilon_{H}}{\lambda_{\min}^{+}(\mathbf{L})} \norm{\mathbf{d}}^{2}\\
        \land~~\bm{1}_{m}^{\top} \otimes \mathbf{I}_{n} \cdot \boldsymbol{\theta}^{k} = \bm{r}\Bigg] \ge 1 - p,
    \end{multline*}
     where $\sqrt{\hat{\mathbf{L}}} = \sqrt{\mathbf{L}}\otimes \mathbf{I}_{n}$ with Laplacian matrix $\mathbf{L} = \sqrt{\mathbf{L}} \cdot \sqrt{\mathbf{L}}$ and tangent space $\mathcal{T} = \{\mathbf{d} \in (\mathbb{R}^{n})^{m}~:~\bm{1}_{m}^{\top} \otimes \mathbf{I}_{n} \cdot \mathbf{d} = \bm{0}\}$. Moreover, if the initial point $\boldsymbol{\theta}^{0}$ is feasible, i.e., $\bm{1}_{m}^{\top} \otimes \mathbf{I}_{n} \cdot \boldsymbol{\theta}^{0} = \bm{r}$, then the feasibility condition $\bm{1}_{m}^{\top} \otimes \mathbf{I}_{n} \cdot \boldsymbol{\theta}^{k} = \bm{r}$ holds for all $k \ge 0$ due to the Laplacian structure of~\eqref{eq: NLGD}, ensuring \emph{anytime feasibility} of the allocation.
\end{theorem}

The theorem shows that \textbf{NLGD} reaches approximate second-order stationarity with high probability while preserving feasibility at every iteration. The result follows by analyzing noisy gradient descent on $\Psi_{\boldsymbol{\theta}^{0}}$ and transferring its first- and second-order properties to the original constrained problem.

\begin{remark}
    The global Lipschitz constants are needed only for offline parameter selection. For fixed admissible $\alpha$, the iteration bound scales as $O(\epsilon_g^{-2}\alpha^{-2})$. Applying $\sqrt{\mathbf{L}}$ exactly may require additional local computation and storage, while polynomial filtering trades this cost for additional communication and approximation error.
\end{remark}

The additional cost introduced by the perturbation mechanism consists of both generating the random vectors and applying the square-root Laplacian. Noise generation requires $\mathcal{O}(n)$ local operations per agent, while the cost of applying $\sqrt{\mathbf{L}}$ depends on its implementation: an exact dense realization generally incurs additional local computation and storage, whereas a polynomial-filter approximation trades dense local mixing for additional communication rounds and approximation error.

\section{Proofs}
\label{sec: Proofs}
\subsection{Proof of Proposition \ref{pro: Same sequence}}
\begin{proof}
    By the definition of $\Psi_{\boldsymbol{\theta}^{0}}$ in \eqref{eq: Auxiliary function},
    \begin{gather}
    \label{eq: Psi gradient}
       \nabla \Psi_{\boldsymbol{\theta}^{0}}(\mathbf{x}^{k}) = \sqrt{\hat{\mathbf{L}}} \cdot \nabla F(\boldsymbol{\theta}^{0} + \sqrt{\hat{\mathbf{L}}} \mathbf{x}^{k}) = \sqrt{\hat{\mathbf{L}}} \cdot \nabla F(\boldsymbol{\theta}^{k}).
    \end{gather}
    Thus, the update in \eqref{eq: LGD with auxiliary function} can be reformulated as
    \begin{gather*}
        \mathbf{x}^{k+1} = \mathbf{x}^{k} - \alpha \nabla \Psi_{\boldsymbol{\theta}^{0}}(\mathbf{x}^{k}) = \mathbf{x}^{k} - \alpha \sqrt{\hat{\mathbf{L}}} \cdot \nabla F(\boldsymbol{\theta}^{0} + \sqrt{\hat{\mathbf{L}}} \mathbf{x}^{k}),\\
        \boldsymbol{\theta}^{k+1} = \boldsymbol{\theta}^{0} + \sqrt{\hat{\mathbf{L}}} \mathbf{x}^{k+1}.
    \end{gather*}
    Then, subtracting $\boldsymbol{\theta}^{k}$ from $\boldsymbol{\theta}^{k+1}$ yields
    \begin{gather*}
        \boldsymbol{\theta}^{k+1} = \boldsymbol{\theta}^{k} - \alpha \hat{\mathbf{L}} \nabla F(\boldsymbol{\theta}^{0} + \sqrt{\hat{\mathbf{L}}}\mathbf{x}^{k}) = \boldsymbol{\theta}^{k} - \alpha \hat{\mathbf{L}} \cdot \nabla  F(\boldsymbol{\theta}^{k}).
    \end{gather*}
    Therefore, the sequence $\{\boldsymbol{\theta}^{k}\}$ generated by \eqref{eq: LGD aggregate form} follows the same sequence as gradient descent applied to $\Psi_{\boldsymbol{\theta}^{0}}$ from the same initial point $\boldsymbol{\theta}^{0} \in (\mathbb{R}^{n})^{m}$ and $\mathbf{x}^{0} = \bm{0} \in (\mathbb{R}^{n})^{m}$ with the same fixed step-size $\alpha > 0$ as claimed.
\end{proof}

\subsection{Proof of Proposition \ref{pro: First order property}}
\begin{lemma}
\label{lem: Sufficient descent}
    Let Assumptions \ref{ass: Lipschitz}, \ref{ass: Coercivity} and \ref{ass: Network} hold. Given initial point $\boldsymbol{\theta}^{0} \in (\mathbb{R}^{n})^{m}$, for any fixed step-size $0 < \alpha < 2/L_{\Psi_{\boldsymbol{\theta}^{0}}}^{g}$ with initial iterate $\mathbf{x}^{0} = \bm{0} \in (\mathbb{R}^{n})^{m}$, the sequence $\{\mathbf{x}^{k}\}$ generated by \eqref{eq: LGD with auxiliary function} satisfies that
    \begin{gather*}
        \Psi_{\boldsymbol{\theta}^{0}}(\mathbf{x}^{k+1}) - \Psi_{\boldsymbol{\theta}^{0}}(\mathbf{x}^{k}) \le (-\frac{1}{\alpha} + \frac{L_{\Psi_{\boldsymbol{\theta}^{0}}}^{g}}{2}) \norm{\mathbf{x}^{k+1} - \mathbf{x}^{k}}^{2} \le 0.
    \end{gather*}
\end{lemma}

\begin{proof}
    By Assumption \ref{ass: Lipschitz}, in view of the update in \eqref{eq: LGD with auxiliary function},  applying Taylor’s theorem yields
    \begin{align}
    \label{eq: Descent lemma}
        \begin{aligned}
            \Psi_{\boldsymbol{\theta}^{0}}(\mathbf{x}^{k+1}) - \Psi_{\boldsymbol{\theta}^{0}}(\mathbf{x}^{k}) & \le \nabla \Psi_{\boldsymbol{\theta}^{0}}(\mathbf{x}^{k})^{\top} (\mathbf{x}^{k+1} - \mathbf{x}^{k}) + \frac{L_{\Psi_{\boldsymbol{\theta}^{0}}}^{g}}{2} \norm{\mathbf{x}^{k+1} - \mathbf{x}^{k}}^{2}\\
            & \le (-\frac{1}{\alpha} + \frac{L_{\Psi_{\boldsymbol{\theta}^{0}}}^{g}}{2}) \norm{\mathbf{x}^{k+1} - \mathbf{x}^{k}}^{2}.
        \end{aligned}
    \end{align}
    Thus, for any fixed $0 < \alpha < 2/L_{\Psi_{\boldsymbol{\theta}^{0}}}^{g}$,
    \begin{gather*}
        \Psi_{\boldsymbol{\theta}^{0}}(\mathbf{x}^{k+1}) - \Psi_{\boldsymbol{\theta}^{0}}(\mathbf{x}^{k}) \le 0
    \end{gather*}
    as claimed.
\end{proof}

\begin{proof}[Proof of Proposition~\ref{pro: First order property}]
    Consider the update in \eqref{eq: LGD with auxiliary function}.
    By Lemma~\ref{lem: Sufficient descent} and
    \eqref{eq: Psi Lipschitz}, for any fixed step-size
    \[
        0 < \alpha
        \le
        \frac{1}
        {\|\sqrt{\mathbf{L}}\|^{2} L_{F}^{g}}
        <
        \frac{2}{L_{\Psi_{\boldsymbol{\theta}^{0}}}^{g}},
    \]
    with initial iterate
    $\mathbf{x}^{0}=\bm{0}\in(\mathbb{R}^{n})^{m}$, we have
    \[
        \Psi_{\boldsymbol{\theta}^{0}}(\mathbf{x}^{k+1})
        -
        \Psi_{\boldsymbol{\theta}^{0}}(\mathbf{x}^{k})
        \le
        -
        \left(
        \frac{1}{\alpha}
        -
        \frac{L_{\Psi_{\boldsymbol{\theta}^{0}}}^{g}}{2}
        \right)
        \norm{\mathbf{x}^{k+1}-\mathbf{x}^{k}}^{2}
        \le 0.
    \]

    By Assumption~\ref{ass: Coercivity}, each $f_i$ is coercive
    and continuous, and hence attains a finite global minimum.
    Let $f_i^\star$ denote its global minimum. Therefore,
    \[
        \Psi_{\boldsymbol{\theta}^{0}}(\mathbf{x})
        =
        F\!\left(
        \boldsymbol{\theta}^{0}
        +
        \sqrt{\hat{\mathbf{L}}}\mathbf{x}
        \right)
        \ge
        \sum_{i=1}^{m} f_i^\star
    \]
    for all $\mathbf{x}$. Summing the descent inequality from $\kappa=0$ to $k-1$ gives
    \[
        \left(
        \frac{1}{\alpha}
        -
        \frac{L_{\Psi_{\boldsymbol{\theta}^{0}}}^{g}}{2}
        \right)
        \sum_{\kappa=0}^{k-1}
        \norm{\mathbf{x}^{\kappa+1}-\mathbf{x}^{\kappa}}^{2}
        \le
        \Psi_{\boldsymbol{\theta}^{0}}(\mathbf{x}^{0})
        -
        \Psi_{\boldsymbol{\theta}^{0}}(\mathbf{x}^{k})
        \le
        \Psi_{\boldsymbol{\theta}^{0}}(\mathbf{x}^{0})
        -
        \sum_{i=1}^{m} f_i^\star .
    \]
    Since
    \[
        \frac{1}{\alpha}
        -
        \frac{L_{\Psi_{\boldsymbol{\theta}^{0}}}^{g}}{2}
        >0,
    \]
    it follows that
    \[
        \sum_{\kappa=0}^{\infty}
        \norm{\mathbf{x}^{\kappa+1}-\mathbf{x}^{\kappa}}^{2}
        <\infty.
    \]
    Consequently,
    \[
        \lim_{k\to\infty}
        \norm{\mathbf{x}^{k+1}-\mathbf{x}^{k}}
        =0.
    \]

    By \eqref{eq: LGD with auxiliary function},
    \[
        \mathbf{x}^{k+1}-\mathbf{x}^{k}
        =
        -\alpha
        \nabla
        \Psi_{\boldsymbol{\theta}^{0}}(\mathbf{x}^{k}),
    \]
    and since $\alpha>0$ is fixed,
    \[
        \lim_{k\to\infty}
        \norm{
        \nabla
        \Psi_{\boldsymbol{\theta}^{0}}(\mathbf{x}^{k})
        }
        =0.
    \]
    Finally, by \eqref{eq: Psi gradient},
    \[
        \lim_{k\to\infty}
        \norm{
        \sqrt{\hat{\mathbf{L}}}
        \cdot \nabla F(\boldsymbol{\theta}^{k})
        }
        =0.
    \]

    Moreover, since
    $(\bm{1}_{m}^{\top}\otimes\mathbf{I}_{n})
    \sqrt{\hat{\mathbf{L}}}=\bm{0}^{\top}$
    and the initial point is feasible,
    \[
        \begin{aligned}
        (\bm{1}_{m}^{\top}\otimes\mathbf{I}_{n})
        \boldsymbol{\theta}^{k}
        &=
        (\bm{1}_{m}^{\top}\otimes\mathbf{I}_{n})
        \left(
        \boldsymbol{\theta}^{0}
        +
        \sqrt{\hat{\mathbf{L}}}\mathbf{x}^{k}
        \right)\\
        &=
        (\bm{1}_{m}^{\top}\otimes\mathbf{I}_{n})
        \boldsymbol{\theta}^{0}
        =
        \bm{r},
        \end{aligned}
    \]
    for all $k\ge0$.

    By Proposition~\ref{pro: Same sequence}, the sequence
    $\{\boldsymbol{\theta}^{k}\}$ generated by
    \eqref{eq: LGD aggregate form} coincides with the sequence
    generated through \eqref{eq: LGD with auxiliary function}
    from the same initial point and with the same fixed step-size
    $\alpha>0$, which concludes the proof.
\end{proof}

\subsection{Proof of Proposition \ref{pro: Second order relations}}
\begin{proof}
    Given
    \[
        \boldsymbol{\theta}
        =
        \boldsymbol{\theta}^{0}
        +
        \sqrt{\hat{\mathbf{L}}}\mathbf{x},
    \]
    from condition i), it follows that
    \begin{gather*}
        \norm{
        \sqrt{\hat{\mathbf{L}}}
        \nabla F(\boldsymbol{\theta})
        }
        =
        \norm{
        \nabla
        \Psi_{\boldsymbol{\theta}^{0}}(\mathbf{x})
        }
        \le \epsilon.
    \end{gather*}
    Moreover, since the initial point
    $\boldsymbol{\theta}^{0}$ is feasible and
    \[
        (\bm{1}_{m}^{\top}\otimes\mathbf{I}_{n})
        \sqrt{\hat{\mathbf{L}}}
        =
        \bm{0}^{\top},
    \]
    it follows that
    \begin{gather*}
        (\bm{1}_{m}^{\top}\otimes\mathbf{I}_{n})
        \boldsymbol{\theta}
        =
        (\bm{1}_{m}^{\top}\otimes\mathbf{I}_{n})
        \left(
        \boldsymbol{\theta}^{0}
        +
        \sqrt{\hat{\mathbf{L}}}\mathbf{x}
        \right) =
        (\bm{1}_{m}^{\top}\otimes\mathbf{I}_{n})
        \boldsymbol{\theta}^{0}
        =
        \bm{r}.
    \end{gather*}
    From condition ii), it follows that for all
    $\mathbf{e}\in(\mathbb{R}^{n})^{m}$,
    \begin{gather*}
        \mathbf{e}^{\top}
        \sqrt{\hat{\mathbf{L}}}
        \nabla^{2}F(\boldsymbol{\theta})
        \sqrt{\hat{\mathbf{L}}}
        \mathbf{e}
        =
        \mathbf{e}^{\top}
        \nabla^{2}
        \Psi_{\boldsymbol{\theta}^{0}}(\mathbf{x})
        \mathbf{e}
        \ge
        -\gamma\norm{\mathbf{e}}^{2}.
    \end{gather*}
    By Assumption~\ref{ass: Network},
    \[
        \operatorname{range}
        \left(\sqrt{\hat{\mathbf{L}}}\right)
        =
        \mathcal{T}
        =
        \left\{
        \mathbf{d}\in(\mathbb{R}^{n})^{m}
        :
        (\bm{1}_{m}^{\top}\otimes\mathbf{I}_{n})
        \mathbf{d}
        =
        \bm{0}
        \right\}.
    \]
    Hence, for every $\mathbf{d}\in\mathcal{T}$, we can choose
    $\mathbf{e}\in\mathcal{T}$ such that
    \[
        \mathbf{d}
        =
        \sqrt{\hat{\mathbf{L}}}\mathbf{e}.
    \]
    Since $\mathbf{e}\in\mathcal{T}$ and the smallest positive
    eigenvalue of $\hat{\mathbf{L}}$ is
    $\lambda_{\min}^{+}(\mathbf{L})$, we have
    \begin{gather*}
        \norm{\mathbf{d}}^{2}
        =
        \mathbf{e}^{\top}
        \hat{\mathbf{L}}
        \mathbf{e}
        \ge
        \lambda_{\min}^{+}(\mathbf{L})
        \norm{\mathbf{e}}^{2}.
    \end{gather*}
    Therefore,
    \begin{gather*}
        \mathbf{d}^{\top}
        \nabla^{2}F(\boldsymbol{\theta})
        \mathbf{d}
        \ge
        -\gamma\norm{\mathbf{e}}^{2} \ge
        -
        \frac{\gamma}
        {\lambda_{\min}^{+}(\mathbf{L})}
        \norm{\mathbf{d}}^{2},
    \end{gather*}
    as claimed.
\end{proof}

\begin{figure*}[t]
\centering
\resizebox{\textwidth}{!}{%
\begin{tikzpicture}[
    node distance=0.4cm and 2.5cm,
    every node/.style={
      draw,
      font=\small,
      minimum width=2.6cm,
      minimum height=0.9cm,
      align=center
    },
    every path/.style={draw, thick, -{Latex[length=3mm]}},
    line width=0.7pt,
    box/.style={draw, font=\sffamily, minimum width=4cm, minimum height=2cm, align=center}
]

\node (l2) {Lemma \ref{lem: Decrease}\\ Gradient Descent};
\node[right=of l2] (l3) {Lemma \ref{lem: Improve or Localize}};
\node[right=of l3] (l4) {Lemma \ref{lem: Localization}};

\node[below=of l3] (l56) {Lemma \ref{lem: Dynamic}~\&~Lemma \ref{lem: Dynamic concentration 1}};
\node[below=of l4] (l7) {Lemma \ref{lem: Dynamic concentration 2}\\ Properties of  Coupling Sequences};
\node[below=of l7] (l8) {Lemma \ref{lem: Decrease for saddle by coupling}\\ Saddle Point Escape};

\node[left=of l8] (p2) {Proposition \ref{pro: Approximate second order}\\ Approximate Second-order Stationarity};

\path (l2) edge (l3);
\path (l3) edge (l4);
\path (l2) edge (p2);
\path (l4) edge (l7);
\path (l56) edge (l7);
\path (l56) edge (l8);
\path (l7) edge (l8);
\path (l8) edge (p2);

\end{tikzpicture}%
}
\caption{Logical structure of the proof of Proposition~\ref{pro: Approximate second order}.}
\label{fig: Proof structure}
\end{figure*}
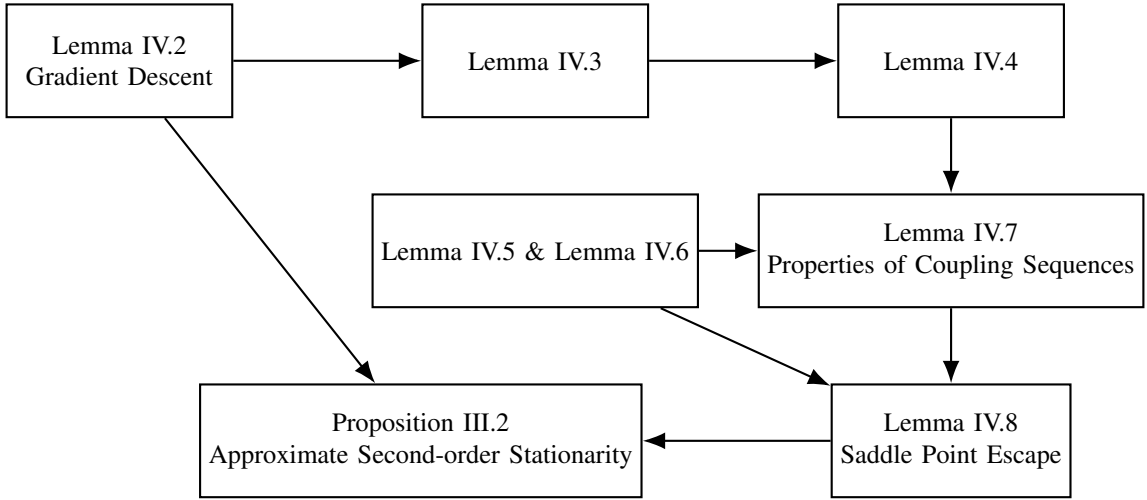

\subsection{Proof of Proposition \ref{pro: Approximate second order}}
Throughout this subsection, $\epsilon_H=\sqrt{\epsilon_gL_{\Psi_{\boldsymbol{\theta}^{0}}}^{H}}$ and $\sigma^{2}=L_{\Psi_{\boldsymbol{\theta}^{0}}}^{g}\alpha\epsilon_g^{2}/(12mn)$ are as specified in Proposition \ref{pro: Approximate second order}.
For the sake of clarity in the forthcoming proofs, we define the following key parameters: given $\rho \ge 1$ (to be specified later), let 
\begin{gather}
    \label{eq: parameters}
    \alpha\coloneqq\frac{1}{L_{\Psi_{\boldsymbol{\theta}^{0}}}^{g}\rho^{7}},
    \qquad
    d\coloneqq\frac{\epsilon_H}{20L_{\Psi_{\boldsymbol{\theta}^{0}}}^{H}\rho^{2}},
    \qquad
    r\coloneqq\left\lceil\frac{\rho}{\alpha\epsilon_H}\right\rceil .
\end{gather}

Before presenting the proof, we introduce two preliminary results used repeatedly below.

\begin{proposition}[Boole--Fr\'echet Inequality]
\label{pro: Boole-frechet Inequality}
    For events $\mathcal{E}_1,\ldots,\mathcal{E}_n$,
    \begin{gather*}
        \max_{1\leq i\leq n}\mathbb{P}(\mathcal{E}_i)\leq\mathbb{P}\left(\bigcup_{i=1}^{n}\mathcal{E}_i\right)\leq\sum_{i=1}^{n}\mathbb{P}(\mathcal{E}_i),
    \end{gather*}
    and consequently
    \begin{gather}
    \label{eq: Boole's inequality}
        \mathbb{P}\left[\bigcap_{i=1}^{n}\mathcal{E}_i\right]\geq1-\sum_{i=1}^{n}\mathbb{P}[\bar{\mathcal{E}}_i]=1-\sum_{i=1}^{n}\left(1-\mathbb{P}[\mathcal{E}_i]\right),
    \end{gather}
    where $\bar{\mathcal{E}}_i$ denotes the complement of $\mathcal{E}_i$.
\end{proposition}

    \begin{proposition}[Material Implication Equivalence, \cite{hurley2011concise}]
    \label{pro: Material implication}
    For any statements $P_A$ and $P_B$, the implication $P_A\Rightarrow P_B$ is logically equivalent to $\bar{P}_A\lor P_B$.
\end{proposition}

For readability, we abuse logical notation when dealing with events. Let $\Omega$ be the whole sample space. Specifically, given events of the form $\mathcal{E}_{i} := \{\omega \in \Omega : P_{i}\} \subseteq \Omega$ and $\mathcal{E}_{j} := \{\omega \in \Omega : P_{j}\} \subseteq \Omega$, where $P_{i}$ and $P_{j}$ are logical predicates defined on $\omega$, we denote:
\begin{gather}
\label{eq: Material implication}
    \mathcal{E}_{i} \Rightarrow \mathcal{E}_{j} := \left\{ \omega \in \Omega :~P_{i} \Rightarrow P_{j} \right\} = \left\{ \omega \in \Omega :~\neg P_{i} \lor P_{j} \right\} = \bar{\mathcal{E}}_{i} \cup \mathcal{E}_{j}.
\end{gather}

The middle step follows the classical logic of material implication (see Proposition~\ref{pro: Material implication}).

The logical structure of the proof of Proposition \ref{pro: Approximate second order} is shown in Fig. \ref{fig: Proof structure}. Lemma \ref{lem: Decrease} bounds the change in $\Psi_{\boldsymbol{\theta}^{0}}$ by separating the gradient-induced decrease from the increase due to random perturbations. Lemma \ref{lem: Improve or Localize} then shows that, with high probability, either the objective decreases sufficiently or the iterates remain localized near the initial point.

\begin{lemma}
\label{lem: Decrease}
    Let Assumptions \ref{ass: Lipschitz} and \ref{ass: Random perturbation} hold. Given $\rho\geq1$, let $\alpha$ and $d$ be defined in \eqref{eq: parameters}. Then, for any $\mathbf{x}^{0}\in(\mathbb{R}^{n})^{m}$ and $t_{0},t\geq0$,
    \begin{multline}
    \label{eq: Decomposition of the change}
    \mathbb{P}\Bigg[\Psi_{\boldsymbol{\theta}^{0}}(\mathbf{x}^{t_{0}+t})-\Psi_{\boldsymbol{\theta}^{0}}(\mathbf{x}^{t_{0}})\leq-\frac{\alpha}{2}\sum_{k=0}^{t-1}\|\nabla\Psi_{\boldsymbol{\theta}^{0}}(\mathbf{x}^{t_{0}+k})\|^{2}\\
    +mn\alpha\sigma^{2}(t+\sqrt{t\rho}+\rho)\Bigg]\geq1-2\mathrm{e}^{-\rho},
    \end{multline}
    where $\{\mathbf{x}^{k}\}$ is generated by \eqref{eq: NLGD with auxiliary function 1}.
\end{lemma}

\begin{proof}
    Since the dynamics in \eqref{eq: NLGD with auxiliary function 1} are time-invariant, it suffices to consider $t_{0}=0$. By Assumption \ref{ass: Lipschitz} and Taylor's theorem, for any $t\geq0$,
    \begin{align*}
        \Psi_{\boldsymbol{\theta}^{0}}(\mathbf{x}^{t+1})-\Psi_{\boldsymbol{\theta}^{0}}(\mathbf{x}^{t})
        &\leq \nabla\Psi_{\boldsymbol{\theta}^{0}}(\mathbf{x}^{t})^{\top}(\mathbf{x}^{t+1}-\mathbf{x}^{t})+\frac{L_{\Psi_{\boldsymbol{\theta}^{0}}}^{g}}{2}\|\mathbf{x}^{t+1}-\mathbf{x}^{t}\|^{2}\\
        &=\left(-\alpha+\frac{L_{\Psi_{\boldsymbol{\theta}^{0}}}^{g}\alpha^{2}}{2}\right)\|\nabla\Psi_{\boldsymbol{\theta}^{0}}(\mathbf{x}^{t})\|^{2}+\frac{L_{\Psi_{\boldsymbol{\theta}^{0}}}^{g}\alpha^{2}}{2}\|\mathbf{n}^{t}\|^{2}\\
        &\quad+\left(-\alpha+L_{\Psi_{\boldsymbol{\theta}^{0}}}^{g}\alpha^{2}\right)\nabla\Psi_{\boldsymbol{\theta}^{0}}(\mathbf{x}^{t})^{\top}\mathbf{n}^{t}.
    \end{align*}
    Since $\rho\geq1$ implies $\alpha=1/(L_{\Psi_{\boldsymbol{\theta}^{0}}}^{g}\rho^{7})\leq1/L_{\Psi_{\boldsymbol{\theta}^{0}}}^{g}$, Cauchy--Schwarz and Young's inequalities give
    \begin{gather*}
        \Psi_{\boldsymbol{\theta}^{0}}(\mathbf{x}^{t+1})-\Psi_{\boldsymbol{\theta}^{0}}(\mathbf{x}^{t})\leq-\frac{\alpha}{2}\|\nabla\Psi_{\boldsymbol{\theta}^{0}}(\mathbf{x}^{t})\|^{2}+\frac{\alpha}{2}\|\mathbf{n}^{t}\|^{2}.
    \end{gather*}
    Summing over $t$ yields
    \begin{gather}
    \label{eq: descent}
        \Psi_{\boldsymbol{\theta}^{0}}(\mathbf{x}^{t})-\Psi_{\boldsymbol{\theta}^{0}}(\mathbf{x}^{0})\leq-\frac{\alpha}{2}\sum_{k=0}^{t-1}\left(\|\nabla\Psi_{\boldsymbol{\theta}^{0}}(\mathbf{x}^{k})\|^{2}-\|\mathbf{n}^{k}\|^{2}\right).
    \end{gather}
    Since $\sum_{k=0}^{t-1}\|\mathbf{n}^{k}/\sigma\|^{2}$ is chi-square distributed with $tmn$ degrees of freedom, Lemma 1 in \cite{laurent2000adaptive} gives
    \begin{gather}
    \label{eq: descent third term}
        \mathbb{P}\left[\sum_{k=0}^{t-1}\|\mathbf{n}^{k}\|^{2}\leq2mn\sigma^{2}(t+\sqrt{t\rho}+\rho)\right]\geq1-2\mathrm{e}^{-\rho}.
    \end{gather}
    Combining \eqref{eq: descent} and \eqref{eq: descent third term} proves \eqref{eq: Decomposition of the change}.
\end{proof}

\begin{lemma}
    \label{lem: Improve or Localize}
    Let Assumptions \ref{ass: Lipschitz} and \ref{ass: Random perturbation} hold. Given $\rho\geq1$, let $\alpha$ and $d$ be defined in \eqref{eq: parameters}. Then, for any $\mathbf{x}^{0}\in(\mathbb{R}^{n})^{m}$ and $t_{0},t\geq0$,
    \begin{multline*}
        \mathbb{P}\Big[\forall\tau\in(0,t],~\|\mathbf{x}^{t_{0}+\tau}-\mathbf{x}^{t_{0}}\|^{2}
        \leq4\alpha\tau\big(\Psi_{\boldsymbol{\theta}^{0}}(\mathbf{x}^{t_{0}})-\Psi_{\boldsymbol{\theta}^{0}}(\mathbf{x}^{t_{0}+\tau})\big)\\
        +8mn\alpha^{2}\sigma^{2}\tau(\tau+\sqrt{\tau\rho}+\rho)\Big]\geq1-4t\mathrm{e}^{-\rho},
    \end{multline*}
    where $\{\mathbf{x}^{k}\}$ is generated by \eqref{eq: NLGD with auxiliary function 1}.
\end{lemma}

\begin{proof}
    By time invariance, it suffices to consider $t_{0}=0$. Define
    \begin{gather*}
        \mathcal{E}_{a,\tau}:=\left\{\sum_{k=0}^{\tau-1}\|\nabla\Psi_{\boldsymbol{\theta}^{0}}(\mathbf{x}^{k})\|^{2}
        \leq2\alpha^{-1}\big(\Psi_{\boldsymbol{\theta}^{0}}(\mathbf{x}^{0})-\Psi_{\boldsymbol{\theta}^{0}}(\mathbf{x}^{\tau})\big)
        +2mn\sigma^{2}(\tau+\sqrt{\tau\rho}+\rho)\right\},
    \end{gather*}
    and
    \begin{gather*}
        \mathcal{E}_{b,\tau}:=\left\{\sum_{k=0}^{\tau-1}\|\mathbf{n}^{k}\|^{2}
        \leq2mn\sigma^{2}(\tau+\sqrt{\tau\rho}+\rho)\right\}.
    \end{gather*}
    By Lemma \ref{lem: Decrease},
    \begin{gather}
        \label{eq: localize first term}
        \mathbb{P}[\mathcal{E}_{a,\tau}]\geq1-2\mathrm{e}^{-\rho}.
    \end{gather}
    Moreover, by Cauchy--Schwarz,
    \begin{align}
        \label{eq: localize}
        \|\mathbf{x}^{\tau}-\mathbf{x}^{0}\|^{2}
        &=\alpha^{2}\left\|\sum_{k=0}^{\tau-1}\left(\nabla\Psi_{\boldsymbol{\theta}^{0}}(\mathbf{x}^{k})+\mathbf{n}^{k}\right)\right\|^{2}\nonumber\\
        &\leq2\alpha^{2}\tau\left(\sum_{k=0}^{\tau-1}\|\nabla\Psi_{\boldsymbol{\theta}^{0}}(\mathbf{x}^{k})\|^{2}
        +\sum_{k=0}^{\tau-1}\|\mathbf{n}^{k}\|^{2}\right).
    \end{align}
    Hence, on $\mathcal{E}_{a,\tau}\cap\mathcal{E}_{b,\tau}$,
    \begin{gather*}
        \|\mathbf{x}^{\tau}-\mathbf{x}^{0}\|^{2}
        \leq4\alpha\tau\big(\Psi_{\boldsymbol{\theta}^{0}}(\mathbf{x}^{0})-\Psi_{\boldsymbol{\theta}^{0}}(\mathbf{x}^{\tau})\big) +8mn\alpha^{2}\sigma^{2}\tau(\tau+\sqrt{\tau\rho}+\rho).
    \end{gather*}
    Therefore, by \eqref{eq: Boole's inequality}, \eqref{eq: descent third term}, and \eqref{eq: localize first term},
    \begin{align*}
        \mathbb{P}\left[\bigcap_{\tau=1}^{t}\left(\mathcal{E}_{a,\tau}\cap\mathcal{E}_{b,\tau}\right)\right]
        &\geq1-\sum_{\tau=1}^{t}\big(1-\mathbb{P}[\mathcal{E}_{a,\tau}]\big)
        -\sum_{\tau=1}^{t}\big(1-\mathbb{P}[\mathcal{E}_{b,\tau}]\big)\\
        &\geq1-4t\mathrm{e}^{-\rho},
    \end{align*}
    which proves the result.
\end{proof}

Let $\gamma^{k}\coloneqq-\lambda_{\min}(\nabla^{2}\Psi_{\boldsymbol{\theta}^{0}}(\mathbf{x}^{k}))$ for $k\geq0$. Whenever $\gamma^{0}>0$, let $\mathbf{e}_{\alpha}$ be a unit eigenvector of $\nabla^{2}\Psi_{\boldsymbol{\theta}^{0}}(\mathbf{x}^{0})$ associated with the eigenvalue $-\gamma^{0}$.

\begin{definition}[Coupling sequences]
    \label{def: Coupling sequence}
    Let $\{(\boldsymbol{\xi}^{k},\zeta^{k})\}_{k\geq0}$ be mutually independent across $k$, with $\boldsymbol{\xi}^{k}$ independent of $\zeta^{k}$ for each $k$, and
    \[
        \boldsymbol{\xi}^{k}\sim\mathcal{N}\!\left(\bm{0},\sigma^{2}\big(\mathbf{I}_{mn}-\mathbf{e}_{\alpha}\mathbf{e}_{\alpha}^{\top}\big)\right),
        \qquad
        \zeta^{k}\sim\mathcal{N}(0,\sigma^{2}).
    \]
    Two sequences $\{\mathbf{y}^{k}\}$ and $\{\mathbf{z}^{k}\}$ generated by \eqref{eq: NLGD with auxiliary function 1} from $\mathbf{y}^{0}=\mathbf{z}^{0}=\mathbf{x}^{0}$ are called coupling sequences if their perturbations are
    \[
        \mathbf{n}_{\mathbf{y}}^{k}=\boldsymbol{\xi}^{k}+\zeta^{k}\mathbf{e}_{\alpha},
        \qquad
        \mathbf{n}_{\mathbf{z}}^{k}=\boldsymbol{\xi}^{k}-\zeta^{k}\mathbf{e}_{\alpha}.
    \]
    Thus, each perturbation sequence has marginal law $\mathcal{N}(\bm{0},\sigma^{2}\mathbf{I}_{mn})$, while $\delta_{\mathbf{n}}^{k}\coloneqq\mathbf{n}_{\mathbf{y}}^{k}-\mathbf{n}_{\mathbf{z}}^{k}=2\zeta^{k}\mathbf{e}_{\alpha}$.
\end{definition}

For the coupling argument below, the current starting point is relabeled as time $0$. For $t\geq0$, define
\begin{gather}
    \label{eq: varsigma}
    \varsigma^{0}(t)\coloneqq\sqrt{4\frac{(1+\alpha\gamma^{0})^{2t}}{2\alpha\gamma^{0}+(\alpha\gamma^{0})^{2}}\sigma^{2}}.
\end{gather}
To analyze coupling sequences near a saddle point, define
\begin{gather}
    \label{eq: ls}
    \ell_{s}\coloneqq\frac{d^{2}}{4\alpha r}-2mn\alpha\sigma^{2}(r+\sqrt{r\rho}+\rho),
\end{gather}
where $\alpha$, $d$, and $r$ are defined in \eqref{eq: parameters}, and
\begin{multline}
    \label{eq: Event A}
    \mathcal{E}_{A}^{0,t}\coloneqq\Big\{\exists\tau\in(0,t]:~\min\big\{\Psi_{\boldsymbol{\theta}^{0}}(\mathbf{y}^{\tau})-\Psi_{\boldsymbol{\theta}^{0}}(\mathbf{y}^{0}),\\
    \Psi_{\boldsymbol{\theta}^{0}}(\mathbf{z}^{\tau})-\Psi_{\boldsymbol{\theta}^{0}}(\mathbf{z}^{0})\big\}\leq-\ell_{s}\Big\},
\end{multline}
\begin{gather}
    \label{eq: Event B}
    \mathcal{E}_{B}^{0,t}\coloneqq\Big\{\forall\tau\in(0,t]:~\max\big\{\|\mathbf{y}^{\tau}-\mathbf{y}^{0}\|,\|\mathbf{z}^{\tau}-\mathbf{z}^{0}\|\big\}\leq d\Big\}.
\end{gather}

The next lemmas establish saddle-point escape using coupling sequences. Lemma \ref{lem: Localization} shows that either one sequence achieves sufficient decrease or both remain localized within $r$ iterations. Lemma \ref{lem: Dynamic} decomposes their difference into noise and Hessian components, analyzed in Lemmas \ref{lem: Dynamic concentration 1} and \ref{lem: Dynamic concentration 2}. These results yield the high-probability escape guarantee in Lemma \ref{lem: Decrease for saddle by coupling}.

\begin{lemma}
    \label{lem: Localization}
    Let Assumptions \ref{ass: Lipschitz} and \ref{ass: Random perturbation} hold. Given $\rho\geq1$, let $\alpha$, $d$, and $r$ be defined in \eqref{eq: parameters}. If $\{\mathbf{y}^{k}\}$ and $\{\mathbf{z}^{k}\}$ are coupling sequences, then
    \begin{gather}
        \label{eq: Localization}
        \mathbb{P}\Big[\mathcal{E}_{A}^{0,r}\cup\mathcal{E}_{B}^{0,r}\Big]\geq1-8r\mathrm{e}^{-\rho}.
    \end{gather}
\end{lemma}

\begin{proof}
Applying Lemma \ref{lem: Improve or Localize} to the coupling sequences gives, for any $t>0$,
    \begin{multline}
        \label{eq: Localize y}
        \mathbb{P}\Big[\forall\tau\in(0,t],~\|\mathbf{y}^{\tau}-\mathbf{y}^{0}\|^{2}
        \leq4\alpha\tau\big(\Psi_{\boldsymbol{\theta}^{0}}(\mathbf{y}^{0})-\Psi_{\boldsymbol{\theta}^{0}}(\mathbf{y}^{\tau})\big)\\
        +8mn\alpha^{2}\sigma^{2}\tau(\tau+\sqrt{\tau\rho}+\rho)\Big]\geq1-4t\mathrm{e}^{-\rho},
    \end{multline}
    and similarly
    \begin{multline}
        \label{eq: Localize z}
        \mathbb{P}\Big[\forall\tau\in(0,t],~\|\mathbf{z}^{\tau}-\mathbf{z}^{0}\|^{2}
        \leq4\alpha\tau\big(\Psi_{\boldsymbol{\theta}^{0}}(\mathbf{z}^{0})-\Psi_{\boldsymbol{\theta}^{0}}(\mathbf{z}^{\tau})\big)\\
        +8mn\alpha^{2}\sigma^{2}\tau(\tau+\sqrt{\tau\rho}+\rho)\Big]\geq1-4t\mathrm{e}^{-\rho}.
    \end{multline}
    By \eqref{eq: Boole's inequality},
    \begin{multline}
        \label{eq: coupling sequences 1}
        \mathbb{P}\Big[\forall\tau\in(0,t],~
        \max\{\|\mathbf{y}^{\tau}-\mathbf{y}^{0}\|^{2},\|\mathbf{z}^{\tau}-\mathbf{z}^{0}\|^{2}\}\leq4\alpha\tau\\
        \Big(\max\{\Psi_{\boldsymbol{\theta}^{0}}(\mathbf{y}^{0})-\Psi_{\boldsymbol{\theta}^{0}}(\mathbf{y}^{\tau}),
        \Psi_{\boldsymbol{\theta}^{0}}(\mathbf{z}^{0})-\Psi_{\boldsymbol{\theta}^{0}}(\mathbf{z}^{\tau})\}
        +2mn\alpha\sigma^{2}(\tau+\sqrt{\tau\rho}+\rho)\Big)\Big]\geq1-8t\mathrm{e}^{-\rho}.
    \end{multline}
    Define
    \begin{multline}
        \label{eq: D1}
        \mathcal{C}^{r}\coloneqq\Big\{\forall\tau\in(0,r],~
        \max\{\|\mathbf{y}^{\tau}-\mathbf{y}^{0}\|^{2},\|\mathbf{z}^{\tau}-\mathbf{z}^{0}\|^{2}\}\\
        \leq4\alpha\tau\Big(\max\{\Psi_{\boldsymbol{\theta}^{0}}(\mathbf{y}^{0})-\Psi_{\boldsymbol{\theta}^{0}}(\mathbf{y}^{\tau}),
        \Psi_{\boldsymbol{\theta}^{0}}(\mathbf{z}^{0})-\Psi_{\boldsymbol{\theta}^{0}}(\mathbf{z}^{\tau})\}
        +2mn\alpha\sigma^{2}(\tau+\sqrt{\tau\rho}+\rho)\Big)\Big\}.
    \end{multline}
    By \eqref{eq: coupling sequences 1}, $\mathbb{P}[\mathcal{C}^{r}]\geq1-8r\mathrm{e}^{-\rho}$. On $\mathcal{C}^{r}\cap\bar{\mathcal{E}}_{A}^{0,r}$, Definition \eqref{eq: Event A} gives, for every $\tau\in(0,r]$,
    \begin{gather*}
        \max\{\Psi_{\boldsymbol{\theta}^{0}}(\mathbf{y}^{0})-\Psi_{\boldsymbol{\theta}^{0}}(\mathbf{y}^{\tau}),
        \Psi_{\boldsymbol{\theta}^{0}}(\mathbf{z}^{0})-\Psi_{\boldsymbol{\theta}^{0}}(\mathbf{z}^{\tau})\}<\ell_{s}.
    \end{gather*}
    Hence, using $\tau\leq r$ and \eqref{eq: ls},
    \begin{align*}
        \max\{\|\mathbf{y}^{\tau}-\mathbf{y}^{0}\|^{2},\|\mathbf{z}^{\tau}-\mathbf{z}^{0}\|^{2}\}
        &\leq4\alpha\tau\Big(\ell_{s}+2mn\alpha\sigma^{2}(\tau+\sqrt{\tau\rho}+\rho)\Big)\\
        &\leq4\alpha r\Big(\ell_{s}+2mn\alpha\sigma^{2}(r+\sqrt{r\rho}+\rho)\Big)=d^{2}.
    \end{align*}
    Thus, $\mathcal{C}^{r}\cap\bar{\mathcal{E}}_{A}^{0,r}\subseteq\mathcal{E}_{B}^{0,r}$. By \eqref{eq: Material implication}, this implies $\mathcal{C}^{r}\subseteq\mathcal{E}_{A}^{0,r}\cup\mathcal{E}_{B}^{0,r}$. Therefore,
    \begin{gather*}
        \mathbb{P}\Big[\mathcal{E}_{A}^{0,r}\cup\mathcal{E}_{B}^{0,r}\Big]\geq\mathbb{P}[\mathcal{C}^{r}]\geq1-8r\mathrm{e}^{-\rho},
    \end{gather*}
    which proves the result.
\end{proof}

Define $\Delta^{k}\coloneqq\mathbf{y}^{k}-\mathbf{z}^{k}$, $\delta_{\mathbf{n}}^{k}\coloneqq\mathbf{n}_{\mathbf{y}}^{k}-\mathbf{n}_{\mathbf{z}}^{k}$,
\begin{gather*}
    \mathcal{I}^{k}\coloneqq\int_{0}^{1}\nabla^{2}\Psi_{\boldsymbol{\theta}^{0}}(s\mathbf{y}^{k}+(1-s)\mathbf{z}^{k})\,\textup{d}s,\qquad \mathcal{H}^{0}\coloneqq\nabla^{2}\Psi_{\boldsymbol{\theta}^{0}}(\mathbf{x}^{0}),
\end{gather*}
and
\begin{gather}
    \label{eq: Dynamic expression}
    \begin{aligned}
        \Delta_{1}^{k}&\coloneqq-\alpha\sum_{\tau=0}^{k-1}(\mathbf{I}_{mn}-\alpha\mathcal{H}^{0})^{k-1-\tau}(\mathcal{I}^{\tau}-\mathcal{H}^{0})\Delta^{\tau},\\
        \Delta_{2}^{k}&\coloneqq-\alpha\sum_{\tau=0}^{k-1}(\mathbf{I}_{mn}-\alpha\mathcal{H}^{0})^{k-1-\tau}\delta_{\mathbf{n}}^{\tau}.
    \end{aligned}
\end{gather}
For $t\geq0$, define
\begin{gather}
    \label{eq: Event C}
    \mathcal{E}_{C}^{0,t}\coloneqq\left\{\forall\tau\in(0,t]:~\|\Delta_{1}^{\tau}\|\leq\frac{\alpha\varsigma^{0}(t)}{10}\right\}.
\end{gather}

\begin{lemma}
    \label{lem: Dynamic}
    Let Assumption \ref{ass: Random perturbation} hold and let $\alpha$ be defined in \eqref{eq: parameters}. If $\{\mathbf{y}^{k}\}$ and $\{\mathbf{z}^{k}\}$ are coupling sequences, then for any $k\geq0$,
    \begin{gather*}
        \Delta^{k}=\Delta_{1}^{k}+\Delta_{2}^{k},
    \end{gather*}
    where $\Delta_{1}^{k}$ and $\Delta_{2}^{k}$ are defined in \eqref{eq: Dynamic expression}.
\end{lemma}

\begin{proof}
    From \eqref{eq: NLGD with auxiliary function 1},
    \begin{align*}
        \Delta^{k+1}
        &=\Delta^{k}-\alpha\left(\nabla\Psi_{\boldsymbol{\theta}^{0}}(\mathbf{y}^{k})-\nabla\Psi_{\boldsymbol{\theta}^{0}}(\mathbf{z}^{k})+\delta_{\mathbf{n}}^{k}\right)\\
        &=\Delta^{k}-\alpha\left(\mathcal{I}^{k}\Delta^{k}+\delta_{\mathbf{n}}^{k}\right)\\
        &=(\mathbf{I}_{mn}-\alpha\mathcal{H}^{0})\Delta^{k}-\alpha\left((\mathcal{I}^{k}-\mathcal{H}^{0})\Delta^{k}+\delta_{\mathbf{n}}^{k}\right).
    \end{align*}
    Hence, for $0\leq\tau\leq k-1$, multiplying by $(\mathbf{I}_{mn}-\alpha\mathcal{H}^{0})^{k-\tau-1}$ gives
    \begin{multline}
        \label{eq: Delat t}
        (\mathbf{I}_{mn}-\alpha\mathcal{H}^{0})^{k-\tau-1}\Delta^{\tau+1}
        =(\mathbf{I}_{mn}-\alpha\mathcal{H}^{0})^{k-\tau}\Delta^{\tau}\\
        -\alpha(\mathbf{I}_{mn}-\alpha\mathcal{H}^{0})^{k-\tau-1}
        \left((\mathcal{I}^{\tau}-\mathcal{H}^{0})\Delta^{\tau}+\delta_{\mathbf{n}}^{\tau}\right).
    \end{multline}
    Summing \eqref{eq: Delat t} over $\tau=0,\ldots,k-1$ and using $\Delta^{0}=0$ yields
    \begin{gather*}
        \Delta^{k}=-\alpha\sum_{\tau=0}^{k-1}(\mathbf{I}_{mn}-\alpha\mathcal{H}^{0})^{k-\tau-1}
        \left((\mathcal{I}^{\tau}-\mathcal{H}^{0})\Delta^{\tau}+\delta_{\mathbf{n}}^{\tau}\right)
        =\Delta_{1}^{k}+\Delta_{2}^{k},
    \end{gather*}
    as claimed.
\end{proof}

\begin{lemma}
    \label{lem: Dynamic concentration 1}
    Let Assumption \ref{ass: Random perturbation} hold and $\epsilon_{H}>0$. There exists $\rho_{\min,1}\geq1$ such that for any $\rho\geq\rho_{\min,1}$, with $\alpha$ and $r$ defined in \eqref{eq: parameters}, the following holds: if $\lambda_{\min}(\nabla^{2}\Psi_{\boldsymbol{\theta}^{0}}(\mathbf{x}^{0}))<-\epsilon_{H}$ and $\{\mathbf{y}^{k}\}$ and $\{\mathbf{z}^{k}\}$ are coupling sequences, then for any $k\geq0$,
    \begin{gather*}
        \mathbb{P}\left[\|\Delta_{2}^{k}\|\leq\alpha\varsigma^{0}(k)\sqrt{2\rho}\right]\geq1-2\mathrm{e}^{-\rho},\\
        \mathbb{P}\left[\|\Delta_{2}^{r}\|\geq\frac{\alpha\varsigma^{0}(r)}{5}\right]\geq\frac{2}{3},
    \end{gather*}
    where $\varsigma^{0}(k)$ is defined in \eqref{eq: varsigma}.
\end{lemma}

\begin{proof}
    By Definition \ref{def: Coupling sequence}, $\delta_{\mathbf{n}}^{k}$ lies in $\mathrm{span}(\mathbf{e}_{\alpha})$. Taking $\mathbf{e}_{\alpha}$ to have unit norm, define
    \begin{gather}
        \label{eq: Sum of dev}
        S_{k}\coloneqq\mathbf{e}_{\alpha}^{\top}\sum_{\tau=0}^{k-1}\left(\mathbf{I}_{mn}-\alpha\nabla^{2}\Psi_{\boldsymbol{\theta}^{0}}(\mathbf{x}^{0})\right)^{k-1-\tau}\delta_{\mathbf{n}}^{\tau}.
    \end{gather}
    Since $\nabla^{2}\Psi_{\boldsymbol{\theta}^{0}}(\mathbf{x}^{0})\mathbf{e}_{\alpha}=-\gamma^{0}\mathbf{e}_{\alpha}$ and the perturbations are Gaussian, $S_{k}$ is Gaussian with zero mean and variance
    \begin{gather*}
        (\Sigma_{k})^{2}=4\frac{(1+\alpha\gamma^{0})^{2k}-1}{2\alpha\gamma^{0}+(\alpha\gamma^{0})^{2}}\sigma^{2}\leq(\varsigma^{0}(k))^{2}.
    \end{gather*}
    Moreover, $\Delta_{2}^{k}=-\alpha S_{k}\mathbf{e}_{\alpha}$, so $\|\Delta_{2}^{k}\|=\alpha|S_{k}|$. The first inequality therefore follows from the Gaussian concentration bound.

    Since $\alpha,\gamma^{0}>0$, there exists $\rho_{\min,1}\geq1$ such that for all $\rho\geq\rho_{\min,1}$,
    \begin{gather*}
        (1+\alpha\gamma^{0})^{2r}-1\geq\frac{1}{3}(1+\alpha\gamma^{0})^{2r},
    \end{gather*}
    and hence $\Sigma_{r}\geq\varsigma^{0}(r)/\sqrt{3}$. Letting $\rho^{\prime}=\sqrt{3}/5$ and using $S_{r}/\Sigma_{r}\sim\mathcal{N}(0,1)$,
    \begin{align*}
        \mathbb{P}\left[\|\Delta_{2}^{r}\|\geq\frac{\alpha\varsigma^{0}(r)}{5}\right]
        &\geq\mathbb{P}\left[|S_{r}|\geq\Sigma_{r}\rho^{\prime}\right]\\
        &\geq1-\frac{2\rho^{\prime}}{\sqrt{2\pi}}
        \geq\frac{2}{3},
    \end{align*}
    which proves the result.
\end{proof}

\begin{lemma}
    \label{lem: Dynamic concentration 2}
    Let Assumptions \ref{ass: Lipschitz} and \ref{ass: Random perturbation} hold and $\epsilon_H>0$. There exists $\rho_{\min,3}\newline\geq1$ such that, for any $\rho\geq\rho_{\min,3}$, with $\alpha$, $d$, and $r$ defined in \eqref{eq: parameters}, the following holds. If $\lambda_{\min}(\nabla^{2}\Psi_{\boldsymbol{\theta}^{0}}(\mathbf{x}^{0}))<-\epsilon_H$ and $\{\mathbf{y}^{k}\}$ and $\{\mathbf{z}^{k}\}$ are coupling sequences, then
    \begin{gather*}
        \mathbb{P}\Big[\mathcal{E}_{A}^{0,r}\cup\mathcal{E}_{C}^{0,r}\Big]\geq1-(2r^{2}+8r)\mathrm{e}^{-\rho}.
    \end{gather*}
\end{lemma}

\begin{proof}
    Lemma \ref{lem: Dynamic concentration 1} and the union bound give, for any $t\geq0$,
    \begin{gather}
        \label{eq: Delta 2 for t times}
        \mathbb{P}\left[\forall\tau\in(0,t],~\|\Delta_{2}^{\tau}\|\leq\alpha\varsigma^{0}(\tau)\sqrt{2\rho}\right]\geq1-2t\mathrm{e}^{-\rho}.
    \end{gather}
    Define the stronger event
    \begin{gather*}
        \mathcal{D}^{t}\coloneqq\left\{\forall\tau\in(0,t],~\|\Delta_{1}^{\tau}\|\leq\frac{\alpha\varsigma^{0}(\tau)}{10}\right\}.
    \end{gather*}
    Since $\varsigma^{0}(\tau)$ is nondecreasing, $\mathcal{D}^{t}\subseteq\mathcal{E}_{C}^{0,t}$. We prove by induction that, for $t\leq r$,
    \begin{gather}
        \label{eq: B to C}
        \mathbb{P}\Big[\mathcal{E}_{B}^{0,t}\Rightarrow\mathcal{D}^{t}\Big]\geq1-2t^{2}\mathrm{e}^{-\rho}.
    \end{gather}
    The claim is immediate for $t=0$ since $\Delta_{1}^{0}=\Delta_{2}^{0}=\Delta^{0}=\bm{0}$. Suppose it holds for $t$, and define
    \begin{gather}
        \label{eq: B to C and}
        \mathcal{G}^{t}\coloneqq\left\{\forall\tau\in(0,t],~\|\Delta_{2}^{\tau}\|\leq\alpha\varsigma^{0}(\tau)\sqrt{2\rho}\right\}.
    \end{gather}
    On $\mathcal{D}^{t}\cap\mathcal{G}^{t}$, Lemma \ref{lem: Dynamic} yields
    \begin{gather*}
        \|\Delta^{\tau}\|\leq\alpha\varsigma^{0}(\tau)\left(\sqrt{2\rho}+\frac{1}{10}\right),\qquad \tau\in(0,t].
    \end{gather*}
    Hence, by \eqref{eq: B to C}, \eqref{eq: Delta 2 for t times}, and the union bound,
    \begin{multline}
        \label{eq: Set B to C and Delta}
        \mathbb{P}\Big[\mathcal{E}_{B}^{0,t}\Rightarrow\big(\mathcal{D}^{t}\cap
        \{\forall\tau\in(0,t]:~\|\Delta^{\tau}\|\leq\alpha\varsigma^{0}(\tau)(\sqrt{2\rho}+\tfrac{1}{10})\}\big)\Big]\\
        \geq1-2(t^{2}+t)\mathrm{e}^{-\rho}.
    \end{multline}
    By Hessian Lipschitz continuity, $\|\mathcal{I}^{\tau}-\mathcal{H}^{0}\|\leq L_{\Psi_{\boldsymbol{\theta}^{0}}}^{H}\max\{\|\mathbf{y}^{\tau}-\mathbf{y}^{0}\|,\|\mathbf{z}^{\tau}-\mathbf{z}^{0}\|\}$. Moreover, $\alpha L_{\Psi_{\boldsymbol{\theta}^{0}}}^{g}\leq1$ implies $\|\mathbf{I}_{mn}-\alpha\mathcal{H}^{0}\|\leq1+\alpha\gamma^{0}$. Hence, by \eqref{eq: Dynamic expression},
    \begin{gather}
        \label{eq: Delta 1}
        \|\Delta_{1}^{t+1}\|\leq\alpha\sum_{\tau=0}^{t}(1+\alpha\gamma^{0})^{t-\tau}L_{\Psi_{\boldsymbol{\theta}^{0}}}^{H}
        \max\{\|\mathbf{y}^{\tau}-\mathbf{y}^{0}\|,\|\mathbf{z}^{\tau}-\mathbf{z}^{0}\|\}\|\Delta^{\tau}\|.
    \end{gather}
    Moreover,
    \begin{gather}
        \label{eq: B to B at t+1}
        \mathcal{E}_{B}^{0,t+1}\subseteq\mathcal{E}_{B}^{0,t}.
    \end{gather}
    Since $\Delta^{0}=\bm{0}$ and $\varsigma^{0}(\tau)$ is proportional to $(1+\alpha\gamma^{0})^{\tau}$,
    \begin{gather}
        \label{eq: Geometric sequence}
        \sum_{\tau=1}^{t}(1+\alpha\gamma^{0})^{t-\tau}\varsigma^{0}(\tau)=t\varsigma^{0}(t).
    \end{gather}
    Therefore, on the event in \eqref{eq: Set B to C and Delta} and $\mathcal{E}_{B}^{0,t+1}$,
    \begin{gather*}
        \|\Delta_{1}^{t+1}\|\leq\alpha^{2}L_{\Psi_{\boldsymbol{\theta}^{0}}}^{H}dt\varsigma^{0}(t)\left(\sqrt{2\rho}+\frac{1}{10}\right).
    \end{gather*}
    Since $r=\lceil\rho/(\alpha\epsilon_H)\rceil\leq2\rho/(\alpha\epsilon_H)$ for sufficiently large $\rho$,
    \[
        L_{\Psi_{\boldsymbol{\theta}^{0}}}^{H}\alpha r d(10\sqrt{2\rho}+1)
        \leq\frac{10\sqrt{2\rho}+1}{10\rho}.
    \]
    Hence there exists $\rho_{\min,2}\geq1$ such that, for all $\rho\geq\rho_{\min,2}$,
    \begin{gather*}
        d\leq\frac{1}{L_{\Psi_{\boldsymbol{\theta}^{0}}}^{H}\alpha r(10\sqrt{2\rho}+1)}.
    \end{gather*}
    Using $t\leq r-1$ and $\varsigma^{0}(t)\leq\varsigma^{0}(t+1)$ therefore gives
    \begin{gather}
        \label{eq: B to C at t+1}
        \mathbb{P}\Big[\mathcal{E}_{B}^{0,t+1}\Rightarrow\mathcal{D}^{t+1}\Big]
        \geq1-2(t^{2}+t)\mathrm{e}^{-\rho}
        \geq1-2(t+1)^{2}\mathrm{e}^{-\rho}.
    \end{gather}
    Thus, \eqref{eq: B to C} holds for all $t\leq r$. Since $\mathcal{D}^{r}\subseteq\mathcal{E}_{C}^{0,r}$,
    \begin{gather*}
        \mathbb{P}\Big[\mathcal{E}_{B}^{0,r}\Rightarrow\mathcal{E}_{C}^{0,r}\Big]\geq1-2r^{2}\mathrm{e}^{-\rho}.
    \end{gather*}
    By Lemma \ref{lem: Localization},
    \begin{gather*}
        \mathbb{P}\Big[\bar{\mathcal{E}}_{A}^{0,r}\Rightarrow\mathcal{E}_{B}^{0,r}\Big]
        =\mathbb{P}\Big[\mathcal{E}_{A}^{0,r}\cup\mathcal{E}_{B}^{0,r}\Big]
        \geq1-8r\mathrm{e}^{-\rho}.
    \end{gather*}
    Since $(\bar{\mathcal{E}}_{A}^{0,r}\Rightarrow\mathcal{E}_{B}^{0,r})\cap(\mathcal{E}_{B}^{0,r}\Rightarrow\mathcal{E}_{C}^{0,r})\subseteq(\bar{\mathcal{E}}_{A}^{0,r}\Rightarrow\mathcal{E}_{C}^{0,r})$, \eqref{eq: Material implication} and the union bound yield
    \begin{gather*}
        \mathbb{P}\Big[\mathcal{E}_{A}^{0,r}\cup\mathcal{E}_{C}^{0,r}\Big]
        \geq1-(2r^{2}+8r)\mathrm{e}^{-\rho}.
    \end{gather*}
    Taking $\rho_{\min,3}\geq\max\{\rho_{\min,1},\rho_{\min,2}\}$ proves the result.
\end{proof}

\begin{lemma}
    \label{lem: Decrease for saddle by coupling}
    Let Assumptions \ref{ass: Lipschitz} and \ref{ass: Random perturbation} hold, let $\epsilon_H>0$, and let the perturbation variance be given by \eqref{eq: Variance}. There exists $\rho_{\min,5}\geq1$ such that, for any $\rho\geq\rho_{\min,5}$ and any starting point $\mathbf{x}^{0}\in(\mathbb{R}^{n})^{m}$ satisfying
    \[
        \lambda_{\min}\!\left(\nabla^{2}\Psi_{\boldsymbol{\theta}^{0}}(\mathbf{x}^{0})\right)<-\epsilon_H,
    \]
    we have
    \begin{gather*}
        \mathbb{P}\Big[\exists\tau\in\{1,\ldots,r\}:~
        \Psi_{\boldsymbol{\theta}^{0}}(\mathbf{x}^{\tau})-\Psi_{\boldsymbol{\theta}^{0}}(\mathbf{x}^{0})
        \leq-\ell_s\Big]
        \geq\frac{1}{3}-(r^{2}+8r)\mathrm{e}^{-\rho},
    \end{gather*}
    where $\ell_s$ is defined in \eqref{eq: ls}.
\end{lemma}

\begin{proof}
    Let $\{\mathbf{y}^{k}\}$ and $\{\mathbf{z}^{k}\}$ be the reflection coupling from Definition \ref{def: Coupling sequence}. By Lemmas \ref{lem: Dynamic concentration 1} and \ref{lem: Dynamic concentration 2}, for sufficiently large $\rho$,
    \begin{gather}
        \label{eq: Delta two}
        \mathbb{P}\Big[\|\Delta_{2}^{r}\|\geq\frac{\alpha\varsigma^{0}(r)}{5}\Big]\geq\frac{2}{3},
    \end{gather}
    and
    \begin{gather}
        \label{eq: Delta one}
        \mathbb{P}\Big[\mathcal{E}_{A}^{0,r}\cup\mathcal{E}_{C}^{0,r}\Big]
        \geq1-(2r^{2}+8r)\mathrm{e}^{-\rho}.
    \end{gather}
    Since $\gamma^{0}>\epsilon_H$ and $r=\lceil\rho/(\alpha\epsilon_H)\rceil$, we have $\alpha\gamma^{0}r\geq\rho$. Moreover, gradient Lipschitz continuity gives $\gamma^{0}\leq L_{\Psi_{\boldsymbol{\theta}^{0}}}^{g}$, so $\alpha\gamma^{0}\leq\rho^{-7}\leq1$. Using $\log(1+u)\geq u/2$ for $u\in[0,1]$,
    \begin{gather*}
        (1+\alpha\gamma^{0})^{r}
        \geq\exp\!\left(\frac{\alpha\gamma^{0}r}{2}\right)
        \geq\mathrm{e}^{\rho/2}.
    \end{gather*}
    Hence, by \eqref{eq: varsigma},
    \begin{gather*}
        \frac{\alpha\varsigma^{0}(r)}{20}
        \geq
        \frac{\sigma\sqrt{\alpha}}{10\sqrt{3L_{\Psi_{\boldsymbol{\theta}^{0}}}^{g}}}\mathrm{e}^{\rho/2}.
    \end{gather*}
    With $d=\epsilon_H/(20L_{\Psi_{\boldsymbol{\theta}^{0}}}^{H}\rho^{2})$ and $\sigma^{2}$ as in \eqref{eq: Variance}, there exists $\rho_{\min,4}\geq1$ such that for all $\rho\geq\rho_{\min,4}$,
    \begin{gather}
        \label{eq: varsigma lower bound}
        \frac{\alpha\varsigma^{0}(r)}{20}\geq d.
    \end{gather}
    Therefore, on $\mathcal{E}_{C}^{0,r}\cap\{\|\Delta_{2}^{r}\|\geq\alpha\varsigma^{0}(r)/5\}$,
    \begin{gather*}
        \max\{\|\mathbf{y}^{r}-\mathbf{y}^{0}\|,\|\mathbf{z}^{r}-\mathbf{z}^{0}\|\}
        \geq\frac12\|\Delta^{r}\|
        \geq\frac12\big(\|\Delta_{2}^{r}\|-\|\Delta_{1}^{r}\|\big)
        \geq\frac{\alpha\varsigma^{0}(r)}{20}
        \geq d,
    \end{gather*}
    and hence
    \begin{gather}
        \label{eq: Delta 2 and C to B}
        \left\{\|\Delta_{2}^{r}\|\geq\frac{\alpha\varsigma^{0}(r)}{5}\right\}
        \cap\mathcal{E}_{C}^{0,r}
        \subseteq\bar{\mathcal{E}}_{B}^{0,r}.
    \end{gather}

    Let
    \[
        \mathcal{G}_{1}\coloneqq\left\{\|\Delta_{2}^{r}\|\geq\frac{\alpha\varsigma^{0}(r)}{5}\right\},\quad
        \mathcal{G}_{2}\coloneqq\mathcal{E}_{A}^{0,r}\cup\mathcal{E}_{C}^{0,r},\quad
        \mathcal{G}_{3}\coloneqq\mathcal{E}_{A}^{0,r}\cup\mathcal{E}_{B}^{0,r}.
    \]
    By \eqref{eq: Delta 2 and C to B}, $\mathcal{G}_{1}\cap\mathcal{G}_{2}\cap\mathcal{G}_{3}\subseteq\mathcal{E}_{A}^{0,r}$. Thus, by \eqref{eq: Delta two}, \eqref{eq: Delta one}, Lemma \ref{lem: Localization}, and the union bound,
    \begin{gather*}
        \mathbb{P}[\mathcal{E}_{A}^{0,r}]
        \geq\frac{2}{3}-(2r^{2}+16r)\mathrm{e}^{-\rho}.
    \end{gather*}
    Let $\mathcal{A}_{y}$ and $\mathcal{A}_{z}$ denote the sufficient-decrease events for $\mathbf{y}^{k}$ and $\mathbf{z}^{k}$ appearing in \eqref{eq: Event A}. Since $\mathcal{E}_{A}^{0,r}=\mathcal{A}_{y}\cup\mathcal{A}_{z}$ and the reflection coupling gives $\{\mathbf{y}^{k}\}$ and $\{\mathbf{z}^{k}\}$ identical marginal laws, $\mathbb{P}[\mathcal{A}_{y}]=\mathbb{P}[\mathcal{A}_{z}]$. Hence,
    \begin{gather*}
        \mathbb{P}[\mathcal{A}_{y}]
        =\mathbb{P}[\mathcal{A}_{z}]
        \geq\frac{1}{3}-(r^{2}+8r)\mathrm{e}^{-\rho}.
    \end{gather*}
    Each marginal process has the same law as \eqref{eq: NLGD with auxiliary function 1}. Taking $\rho_{\min,5}\geq\max\{\rho_{\min,3},\rho_{\min,4}\}$ proves the result.
\end{proof}

Let $\mathcal{H}_{k}\coloneqq\sigma(\mathbf{n}^{0},\ldots,\mathbf{n}^{k-1})$. The proof of Lemma \ref{lem: Decrease for saddle by coupling} is uniform in its starting point and uses only the future perturbations. Because the recursion is time-homogeneous and $\{\mathbf{n}^{k}\}$ is independent over time, the same bound therefore holds conditionally at any bounded $\mathcal{H}_{k}$-stopping time $T$. In particular, on the event
\[
    \lambda_{\min}\!\left(\nabla^{2}\Psi_{\boldsymbol{\theta}^{0}}(\mathbf{x}^{T})\right)<-\epsilon_H,
\]
we have almost surely
\begin{gather}
    \label{eq: conditional saddle decrease}
    \mathbb{P}\left[\exists\tau\in\{1,\ldots,r\}:~
    \Psi_{\boldsymbol{\theta}^{0}}(\mathbf{x}^{T+\tau})-\Psi_{\boldsymbol{\theta}^{0}}(\mathbf{x}^{T})
    \leq-\ell_s\,\middle|\,\mathcal{H}_{T}\right]\geq q_{\rho},
\end{gather}
where $q_{\rho}\coloneqq\frac{1}{3}-(r^{2}+8r)\mathrm{e}^{-\rho}$.

Finally, we are prepared to prove Proposition \ref{pro: Approximate second order}, which demonstrates that by selecting a sufficiently large confidence parameter $\rho > 0$, after a sufficient number of iterations, at least one of the updates is an approximate second-order stationary point with high probability.

\begin{proof}[Proof of Proposition \ref{pro: Approximate second order}]
    Let
    \[
        \underline{\Psi}\coloneqq\sum_{i=1}^{m}f_i^\star,
        \qquad
        D\coloneqq\Psi_{\boldsymbol{\theta}^{0}}(\mathbf{x}^{0})-\underline{\Psi}\geq0.
    \]
    If $D=0$, then $\mathbf{x}^{0}$ is a global minimizer of the smooth function $\Psi_{\boldsymbol{\theta}^{0}}$, so $\nabla\Psi_{\boldsymbol{\theta}^{0}}(\mathbf{x}^{0})=\bm{0}$ and $\nabla^{2}\Psi_{\boldsymbol{\theta}^{0}}(\mathbf{x}^{0})\succeq0$; the claim is immediate. Hence assume $D>0$.

    Given $\epsilon_g>0$ and $\rho\geq1$, let $\alpha$ be defined in \eqref{eq: parameters} and set
    \begin{gather}
        \label{eq: K}
        K\coloneqq
        \left\lceil
        \frac{D}{L_{\Psi_{\boldsymbol{\theta}^{0}}}^{g}\epsilon_g^{2}\alpha^{2}}
        \right\rceil
        =
        \left\lceil
        D\,\epsilon_g^{-2}\alpha^{-1}\rho^{7}
        \right\rceil .
    \end{gather}
    Define, for $k\geq0$,
    \begin{gather*}
        \mathcal{E}_{1}^{k}\coloneqq
        \{\|\nabla\Psi_{\boldsymbol{\theta}^{0}}(\mathbf{x}^{k})\|\geq\epsilon_g\},\\
        \mathcal{E}_{2}^{k}\coloneqq
        \{\|\nabla\Psi_{\boldsymbol{\theta}^{0}}(\mathbf{x}^{k})\|<\epsilon_g,\ 
        \lambda_{\min}(\nabla^{2}\Psi_{\boldsymbol{\theta}^{0}}(\mathbf{x}^{k}))<-\epsilon_H\},\\
        \mathcal{E}_{3}^{k}\coloneqq
        \{\|\nabla\Psi_{\boldsymbol{\theta}^{0}}(\mathbf{x}^{k})\|<\epsilon_g,\ 
        \lambda_{\min}(\nabla^{2}\Psi_{\boldsymbol{\theta}^{0}}(\mathbf{x}^{k}))\geq-\epsilon_H\}.
    \end{gather*}
    These events form a partition for each $k$. It suffices to show that $\mathcal{E}_{3}^{k}$ occurs at least once for $k=0,\ldots,K-1$ with probability at least $1-p$. Let
    \begin{gather}
        \label{eq: Event P}
        \mathcal{P}\coloneqq
        \left\{\sum_{k=0}^{K-1}
        \left(\mathbf{1}_{\mathcal{E}_{1}^{k}}+\mathbf{1}_{\mathcal{E}_{2}^{k}}\right)<K\right\}
        =
        \left\{\sum_{k=0}^{K-1}\mathbf{1}_{\mathcal{E}_{3}^{k}}\geq1\right\},
    \end{gather}
    and
    \begin{gather*}
        \mathcal{P}_{1}\coloneqq
        \left\{
        \sum_{k=0}^{K-1}\left(\mathbf{1}_{\mathcal{E}_{1}^{k}}+\mathbf{1}_{\mathcal{E}_{2}^{k}}\right)<K
        \ \lor\
        \sum_{k=0}^{K-1}\mathbf{1}_{\mathcal{E}_{1}^{k}}<\frac{K}{2}
        \right\},\\
        \mathcal{P}_{2}\coloneqq
        \left\{
        \sum_{k=0}^{K-1}\left(\mathbf{1}_{\mathcal{E}_{1}^{k}}+\mathbf{1}_{\mathcal{E}_{2}^{k}}\right)<K
        \ \lor\
        \sum_{k=0}^{K-1}\mathbf{1}_{\mathcal{E}_{2}^{k}}<\frac{K}{2}
        \right\}.
    \end{gather*}
    Then
    \begin{gather}
        \label{eq: P1 and P2 to P}
        \mathcal{P}_{1}\cap\mathcal{P}_{2}\subseteq\mathcal{P}.
    \end{gather}

    Define the noise-concentration event
    \begin{gather*}
        \mathcal{N}_{K}\coloneqq
        \left\{\sum_{k=0}^{K-1}\|\mathbf{n}^{k}\|^{2}
        \leq2mn\sigma^{2}(K+\sqrt{K\rho}+\rho)\right\}.
    \end{gather*}
    By \eqref{eq: descent third term},
    \begin{gather}
        \label{eq: global noise event}
        \mathbb{P}[\mathcal{N}_{K}]\geq1-2\mathrm{e}^{-\rho}.
    \end{gather}

    We first bound $\mathbb{P}[\mathcal{P}_{1}]$. On $\bar{\mathcal{P}}_{1}$, at least $K/2$ of the indices $k=0,\ldots,K-1$ satisfy $\mathcal{E}_{1}^{k}$. Hence, by \eqref{eq: descent}, on $\bar{\mathcal{P}}_{1}\cap\mathcal{N}_{K}$,
    \begin{gather*}
        \Psi_{\boldsymbol{\theta}^{0}}(\mathbf{x}^{K})
        -\Psi_{\boldsymbol{\theta}^{0}}(\mathbf{x}^{0})
        \leq-L_{1},
    \end{gather*}
    where
    \begin{gather*}
        L_{1}\coloneqq
        \frac{\alpha K\epsilon_g^{2}}{4}
        -mn\alpha\sigma^{2}(K+\sqrt{K\rho}+\rho).
    \end{gather*}
    Since $D>0$, for sufficiently large $\rho$ we have $K\geq\rho$. Using \eqref{eq: Variance} and $L_{\Psi_{\boldsymbol{\theta}^{0}}}^{g}\alpha=\rho^{-7}$,
    \begin{align*}
        L_{1}
        &\geq
        \alpha K\left(\frac{\epsilon_g^{2}}{4}-3mn\sigma^{2}\right)
        =
        \frac{\alpha K\epsilon_g^{2}}{4}(1-\rho^{-7})
        \geq\frac{\alpha K\epsilon_g^{2}}{8}.
    \end{align*}
    Moreover, by \eqref{eq: K},
    \begin{gather*}
        \frac{\alpha K\epsilon_g^{2}}{8}
        \geq\frac{D}{8L_{\Psi_{\boldsymbol{\theta}^{0}}}^{g}\alpha}
        =\frac{D\rho^{7}}{8}.
    \end{gather*}
    Thus, there exists $\rho_{\min,6}\geq1$ such that $L_{1}>D$ for all $\rho\geq\rho_{\min,6}$. Since $\Psi_{\boldsymbol{\theta}^{0}}(\mathbf{x})\geq\underline{\Psi}$ for all $\mathbf{x}$, the event $\bar{\mathcal{P}}_{1}\cap\mathcal{N}_{K}$ is empty. Therefore,
    \begin{gather}
        \label{eq: Contradiction 1}
        \mathbb{P}[\mathcal{P}_{1}]
        \geq\mathbb{P}[\mathcal{N}_{K}]
        \geq1-2\mathrm{e}^{-\rho}.
    \end{gather}

    We next bound $\mathbb{P}[\mathcal{P}_{2}]$. On $\bar{\mathcal{P}}_{2}$, at least $K/2$ indices in $\{0,\ldots,K-1\}$ satisfy $\mathcal{E}_{2}^{k}$. For sufficiently large $\rho$, $r\geq\rho$ and
    \begin{gather*}
        \ell_s
        \geq
        \frac{d^{2}}{4\alpha r}-6mn\alpha\sigma^{2}r.
    \end{gather*}
    Also, for sufficiently large $\rho$,
    \[
        r\leq\frac{2\rho}{\alpha\epsilon_H}.
    \]
    Using \eqref{eq: Variance}, \eqref{eq: parameters}, and $\epsilon_H^{2}=\epsilon_gL_{\Psi_{\boldsymbol{\theta}^{0}}}^{H}$,
    \begin{gather*}
        \frac{6mn\alpha\sigma^{2}r}{d^{2}/(8\alpha r)}
        =
        \frac{48mn\alpha^{2}\sigma^{2}r^{2}}{d^{2}}
        \leq\frac{6400}{\rho}.
    \end{gather*}
    Hence, there exists $\rho_{\min,7}\geq1$ such that, for all $\rho\geq\rho_{\min,7}$,
    \begin{gather}
        \label{eq: ls lower bound}
        \ell_s\geq\frac{d^{2}}{8\alpha r}>0.
    \end{gather}
    Consequently,
    \begin{gather}
        \label{eq: r/ls}
        \frac{r}{\ell_s}
        \leq\frac{8\alpha r^{2}}{d^{2}}
        \leq12800L_{\Psi_{\boldsymbol{\theta}^{0}}}^{g}\epsilon_g^{-2}\rho^{13}.
    \end{gather}

    Recall
    \[
        q_{\rho}=\frac13-(r^{2}+8r)\mathrm{e}^{-\rho}.
    \]
    Since $r$ is polynomial in $\rho$ and $K/r$ grows as $\rho^{6}$ when $D>0$, there exists $\rho_{\min,8}\geq\max\{\rho_{\min,5},\rho_{\min,7}\}$ such that, for all $\rho\geq\rho_{\min,8}$,
    \begin{gather}
        \label{eq: large rho block conditions}
        q_{\rho}\geq\frac14,\qquad K\geq8r,\qquad K\geq\rho.
    \end{gather}

    We now construct disjoint saddle blocks. Let $T_{1}$ be the first $k\in\{0,\ldots,K-r\}$ for which $\mathcal{E}_{2}^{k}$ occurs; if no such $k$ exists, set $T_{1}=K$. Recursively, if $T_j=K$, set $Z_j=1$, $\tau_j=0$, and $T_{j+1}=K$. Otherwise define
    \begin{gather*}
        \mathcal{A}_{j}\coloneqq
        \left\{\exists\tau\in\{1,\ldots,r\}:~
        \Psi_{\boldsymbol{\theta}^{0}}(\mathbf{x}^{T_j+\tau})
        -\Psi_{\boldsymbol{\theta}^{0}}(\mathbf{x}^{T_j})
        \leq-\ell_s\right\},
    \end{gather*}
    \begin{gather*}
        \tau_j\coloneqq
        \begin{cases}
            \min\{\tau\in\{1,\ldots,r\}:
            \Psi_{\boldsymbol{\theta}^{0}}(\mathbf{x}^{T_j+\tau})
            -\Psi_{\boldsymbol{\theta}^{0}}(\mathbf{x}^{T_j})\leq-\ell_s\},
            & \mathcal{A}_j\text{ occurs},\\
            r, & \text{otherwise},
        \end{cases}
    \end{gather*}
    where
    \begin{gather*}
        Z_j\coloneqq\mathbf{1}_{\mathcal{A}_j}.
    \end{gather*}
    and let $T_{j+1}$ be the first $k\in\{T_j+\tau_j,\ldots,K-r\}$ satisfying $\mathcal{E}_{2}^{k}$; if no such index exists, set $T_{j+1}=K$.

    Under $\bar{\mathcal{P}}_{2}$, at most $r$ saddle indices lie in the final $r$ positions, and each selected block contains at most $r$ indices. Therefore the first
    \begin{gather}
        \label{eq: saddle block number}
        M\coloneqq\left\lfloor\frac{K/2-r}{r}\right\rfloor
    \end{gather}
    selected blocks are genuine, and \eqref{eq: large rho block conditions} gives
    \begin{gather*}
        M\geq\frac{K}{4r}.
    \end{gather*}

    By \eqref{eq: conditional saddle decrease}, for a genuine block,
    \[
        \mathbb{E}[Z_j\mid\mathcal{H}_{T_j}]\geq q_{\rho},
    \]
    while for a dummy block $Z_j=1$, so the same bound holds. Set $\mu_j\coloneqq\mathbb{E}[Z_j\mid\mathcal{H}_{T_j}]$ and define the block filtration $\mathcal{B}_{0}\coloneqq\mathcal{H}_{T_1}$ and $\mathcal{B}_{j}\coloneqq\mathcal{H}_{T_{j+1}}$ for $j\geq1$. Then $Z_j$ is $\mathcal{B}_{j}$-measurable, $\mu_j$ is $\mathcal{B}_{j-1}$-measurable, and
    \[
        \mathbb{E}[Z_j-\mu_j\mid\mathcal{B}_{j-1}]=0.
    \]
    Thus $\{Z_j-\mu_j\}_{j=1}^{M}$ is a martingale-difference sequence with increments bounded in absolute value by $1$. Therefore, by the Azuma--Hoeffding inequality,
    \begin{gather}
        \label{eq: saddle block concentration}
        \mathbb{P}\left[\sum_{j=1}^{M}Z_j<\frac{q_{\rho}M}{2}\right]
        \leq
        \exp\left(-\frac{q_{\rho}^{2}M}{8}\right).
    \end{gather}
    Since $M\geq K/(4r)$ and $K/r$ grows as $\rho^{6}$, there exists $\rho_{\min,9}\geq\rho_{\min,8}$ such that
    \begin{gather}
        \label{eq: saddle block tail}
        \exp\left(-\frac{q_{\rho}^{2}M}{8}\right)\leq\mathrm{e}^{-\rho}
    \end{gather}
    for all $\rho\geq\rho_{\min,9}$.

    On $\bar{\mathcal{P}}_{2}$, all first $M$ blocks are genuine. On the complement of the event in \eqref{eq: saddle block concentration},
    \[
        S\coloneqq\sum_{j=1}^{M}Z_j
        \geq\frac{q_{\rho}M}{2}
        \geq\frac{K}{32r}.
    \]
    The successful block intervals are disjoint, and each contributes an objective decrease of at least $\ell_s$. For every remaining individual iteration, the one-step estimate in the proof of Lemma \ref{lem: Decrease} gives
    \[
        \Psi_{\boldsymbol{\theta}^{0}}(\mathbf{x}^{k+1})
        -\Psi_{\boldsymbol{\theta}^{0}}(\mathbf{x}^{k})
        \leq\frac{\alpha}{2}\|\mathbf{n}^{k}\|^{2}.
    \]
    Hence, on $\bar{\mathcal{P}}_{2}\cap\mathcal{N}_{K}$ and the complement of the event in \eqref{eq: saddle block concentration},
    \begin{align}
        \label{eq: P2 total decrease}
        \Psi_{\boldsymbol{\theta}^{0}}(\mathbf{x}^{K})
        -\Psi_{\boldsymbol{\theta}^{0}}(\mathbf{x}^{0})
        &\leq
        -S\ell_s+\frac{\alpha}{2}\sum_{k=0}^{K-1}\|\mathbf{n}^{k}\|^{2}\nonumber\\
        &\leq
        -\frac{K\ell_s}{32r}
        +mn\alpha\sigma^{2}(K+\sqrt{K\rho}+\rho).
    \end{align}

    The first term dominates the initial objective gap for sufficiently large $\rho$. Indeed, by \eqref{eq: K} and \eqref{eq: r/ls},
    \begin{gather*}
        \frac{K\ell_s}{32r}
        \geq
        \frac{D\rho}{409600}.
    \end{gather*}
    Moreover, since $K\geq\rho$ and $K\leq D/(L_{\Psi_{\boldsymbol{\theta}^{0}}}^{g}\epsilon_g^{2}\alpha^{2})+1$,
    \begin{align*}
        mn\alpha\sigma^{2}(K+\sqrt{K\rho}+\rho)
        &\leq3mn\alpha\sigma^{2}K\\
        &=\frac{L_{\Psi_{\boldsymbol{\theta}^{0}}}^{g}\alpha^{2}\epsilon_g^{2}}{4}K\\
        &\leq\frac{D}{4}
        +\frac{L_{\Psi_{\boldsymbol{\theta}^{0}}}^{g}\alpha^{2}\epsilon_g^{2}}{4}.
    \end{align*}
    Since $D>0$, by increasing $\rho_{\min,9}$ if necessary, the last term is at most $D/4$, and
    \begin{gather*}
        \frac{K\ell_s}{32r}
        >
        D+mn\alpha\sigma^{2}(K+\sqrt{K\rho}+\rho).
    \end{gather*}
    Thus \eqref{eq: P2 total decrease} would imply
    \[
        \Psi_{\boldsymbol{\theta}^{0}}(\mathbf{x}^{K})<\underline{\Psi},
    \]
    contradicting the definition of $\underline{\Psi}$. Hence
    \[
        \bar{\mathcal{P}}_{2}
        \subseteq
        \mathcal{N}_{K}^{c}
        \cup
        \left\{\sum_{j=1}^{M}Z_j<\frac{q_{\rho}M}{2}\right\}.
    \]
    Using \eqref{eq: global noise event}, \eqref{eq: saddle block tail}, and the union bound,
    \begin{gather}
        \label{eq: Contradiction 2}
        \mathbb{P}[\mathcal{P}_{2}]
        \geq1-3\mathrm{e}^{-\rho}.
    \end{gather}

    Combining \eqref{eq: Contradiction 1}, \eqref{eq: Contradiction 2}, and \eqref{eq: P1 and P2 to P} gives
    \begin{gather*}
        \mathbb{P}[\mathcal{P}]
        \geq1-5\mathrm{e}^{-\rho}.
    \end{gather*}
    Hence there exists $\rho_{\min,10}\geq\rho_{\min,9}$ such that for all $\rho\geq\rho_{\min,10}$,
    \begin{multline}
        \label{eq: Approximately second order staionarity}
        \mathbb{P}\Big[\exists k\in[0,K]:~
        \|\nabla\Psi_{\boldsymbol{\theta}^{0}}(\mathbf{x}^{k})\|\leq\epsilon_g
        \ \land\
        \lambda_{\min}(\nabla^{2}\Psi_{\boldsymbol{\theta}^{0}}(\mathbf{x}^{k}))\geq-\epsilon_H\Big]\\
        \geq1-5\mathrm{e}^{-\rho}
        \geq1-\mathrm{e}^{-\rho/2}.
    \end{multline}

    Finally, \eqref{eq: Approximate second order} follows by choosing
    \begin{gather*}
        \bar{\alpha}
        =
        \min\left\{
        \frac{(\rho_{\min,11})^{-7}}{L_{\Psi_{\boldsymbol{\theta}^{0}}}^{g}},
        \frac{2\ln(1/p)}{L_{\Psi_{\boldsymbol{\theta}^{0}}}^{g}},
        \frac{1}{L_{\Psi_{\boldsymbol{\theta}^{0}}}^{g}[2\ln(1/p)]^{7}}
        \right\},
    \end{gather*}
    where
    \begin{gather*}
        \rho_{\min,11}
        \geq
        \max\{\rho_{\min,6},\rho_{\min,8},\rho_{\min,9},\rho_{\min,10}\}.
    \end{gather*}
    For any $0<\alpha\leq\bar{\alpha}$, let
    \[
        \rho
        =
        \left(\frac{1}{L_{\Psi_{\boldsymbol{\theta}^{0}}}^{g}\alpha}\right)^{1/7}.
    \]
    Then $\rho\geq\rho_{\min,11}$ and $\rho\geq2\ln(1/p)$. Hence
    \[
        1-\mathrm{e}^{-\rho/2}\geq1-p,
    \]
    which proves \eqref{eq: Approximate second order}.
\end{proof}

\subsection{Proof of Theorem \ref{the: Approximate second order}}
\begin{proof}
    By \eqref{eq: Approximate second order} in
    Proposition~\ref{pro: Approximate second order}, there exists
    \begin{gather*}
        \bar{\alpha}
        \le
        \min\left\{
        \frac{1}{L_{\Psi_{\boldsymbol{\theta}^{0}}}^{g}},
        -\frac{2\ln(p)}{L_{\Psi_{\boldsymbol{\theta}^{0}}}^{g}}
        \right\}
    \end{gather*}
    such that for any step-size $\alpha \le \bar{\alpha}$,
    with $K$ as per \eqref{eq: Iteration}, and initial condition
    satisfying $\mathbf{x}^{0}=\bm{0}$, it follows that after
    $K$ iterations of \eqref{eq: NLGD with auxiliary function},
    \begin{align*}
        \mathbb{P}\Bigg[
        \exists k\in[0,K],~
        \|\nabla\Psi_{\boldsymbol{\theta}^{0}}(\mathbf{x}^{k})\|
        \le\epsilon_{g}
        ~\land~
        \lambda_{\min}
        \left(
        \nabla^{2}\Psi_{\boldsymbol{\theta}^{0}}(\mathbf{x}^{k})
        \right)
        \ge-\epsilon_{H}
        \Bigg]
        \ge 1-p.
    \end{align*}

    By Proposition~\ref{pro: NGD interpretation},
    \[
        \boldsymbol{\theta}^{k}
        =
        \boldsymbol{\theta}^{0}
        +
        \sqrt{\hat{\mathbf{L}}}\mathbf{x}^{k}.
    \]
    Since $\boldsymbol{\theta}^{0}$ is feasible and
    \[
        (\bm{1}_{m}^{\top}\otimes\mathbf{I}_{n})
        \sqrt{\hat{\mathbf{L}}}
        =
        \bm{0}^{\top},
    \]
    it follows that, for every $k\ge0$,
    \begin{gather*}
        (\bm{1}_{m}^{\top}\otimes\mathbf{I}_{n})
        \boldsymbol{\theta}^{k} =
        (\bm{1}_{m}^{\top}\otimes\mathbf{I}_{n})
        \left(
        \boldsymbol{\theta}^{0}
        +
        \sqrt{\hat{\mathbf{L}}}\mathbf{x}^{k}
        \right) =
        (\bm{1}_{m}^{\top}\otimes\mathbf{I}_{n})
        \boldsymbol{\theta}^{0} =
        \bm{r}.
    \end{gather*}

    Then, by Proposition~\ref{pro: Second order relations}
    and Proposition~\ref{pro: NGD interpretation},
    after $K$ iterations of~\eqref{eq: NLGD},
    \begin{multline*}
        \mathbb{P}\Bigg[
        \exists k\in[0,K],~
        \|\sqrt{\hat{\mathbf{L}}}\cdot
        \nabla F(\boldsymbol{\theta}^{k})\|
        \le\epsilon_{g}\\
        \land~
        \forall\,\mathbf{d}\in\mathcal{T},~
        \mathbf{d}^{\top}
        \nabla^{2}F(\boldsymbol{\theta}^{k})
        \mathbf{d}
        \ge
        -
        \frac{\epsilon_{H}}
        {\lambda_{\min}^{+}(\mathbf{L})}
        \|\mathbf{d}\|^{2}
        ~\land~
        (\bm{1}_{m}^{\top}\otimes\mathbf{I}_{n})
        \boldsymbol{\theta}^{k}
        =
        \bm{r}
        \Bigg]
        \ge 1-p,
    \end{multline*}
    as claimed.
\end{proof}

\section{Numerical Examples}
\label{sec: Numerical Examples}
\subsection{Smart Grid: Load Control and Demand Response}

One motivating example arises in load control and demand response in smart grids. In this setting, each agent acts as a \emph{prosumer}, capable of both consuming and supplying electricity to the grid. The decision variable reflects the net power exchange within a given time window: positive values represent net consumption, while negative values correspond to power generation or injection into the grid. Each agent aims to balance individual benefit with system-wide constraints. The objective captures a trade-off between \emph{diminishing marginal returns} and \emph{increasing marginal cost}, which penalizes excessive net power flow in either direction. This leads to a non-convex distributed resource allocation problem of the form
\begin{gather*}
    \min_{\boldsymbol{\theta} \in (\mathbb{R}^{n})^{m}} \quad
    F(\boldsymbol{\theta}) \triangleq \sum_{i=1}^{m} f_{i}(\boldsymbol{\theta}_{i})
    \quad \text{subject to} \quad
    \sum_{i=1}^{m} \boldsymbol{\theta}_{i} = \mathbf{r},
\end{gather*}
where $\boldsymbol{\theta} = [(\boldsymbol{\theta}_{1})^{\top}, \dots, (\boldsymbol{\theta}_{m})^{\top}]^{\top} \in (\mathbb{R}^{n})^{m}$ is the vector of local decisions and each local cost function is given by
\begin{gather*}
    f_{i}(\boldsymbol{\theta}_{i}) = a_{i} \boldsymbol{\theta}_{i}^{2} - b_{i} \log(1 + \boldsymbol{\theta}_{i}^{2}),
\end{gather*}
with agent-specific parameters $a_{i}, b_{i} > 0$. For admissible choices of $(a_i,b_i)$, the above objective is non-convex and admits strict saddle points, making it suitable for evaluating second-order behavior. The quadratic term models increasing marginal cost, while the logarithmic term introduces diminishing returns. The coupling constraint enforces system-wide power balance.

The smart grid simulation is conducted over a connected Watts--Strogatz communication network \cite{watts1998collective} with parameters $m = 100$,
$n = 1$, neighborhood size $k = 4$, and rewiring probability $p = 0.2$. The stepsize is fixed at $\alpha = 0.001$, and the noise standard deviation is set to $\sigma = 0.05$.
All runs are initialized in a neighborhood of the saddle point $\boldsymbol{\theta}=\mathbf{0}$ to explicitly examine saddle-point escape behavior.

\begin{figure*}
    \centering
    \subfloat[Network.\label{fig: Network}]{%
        \includegraphics[width=0.32\textwidth]{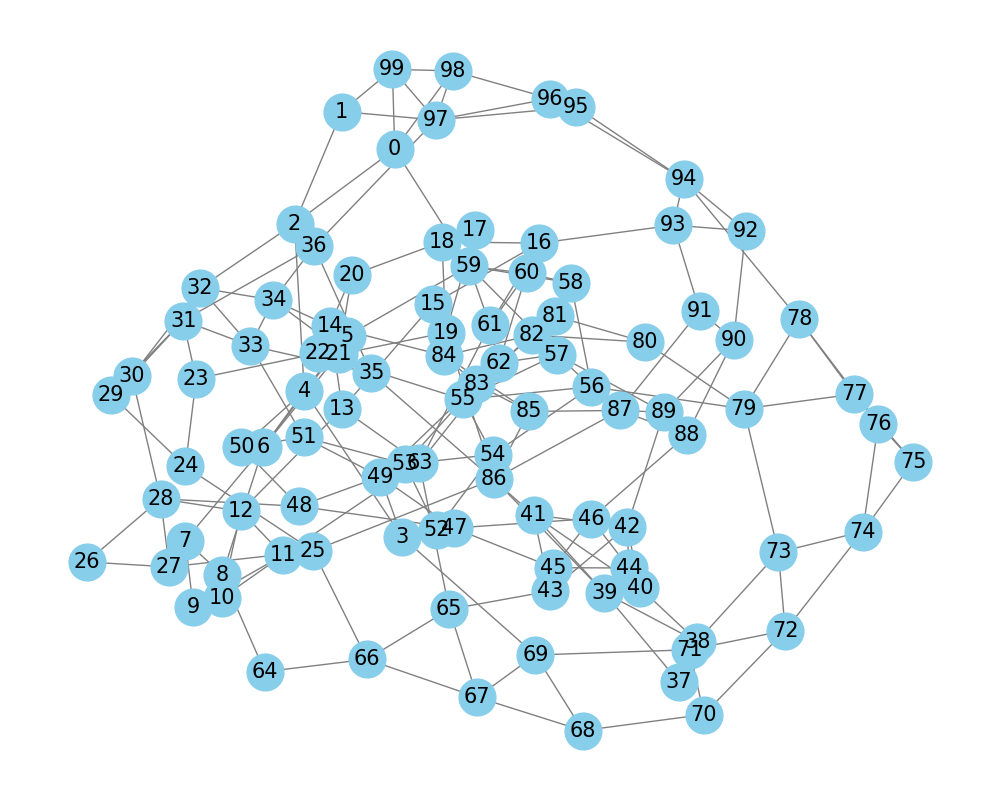}}
    \hfill
    \subfloat[$f(\boldsymbol{\theta}^{k})$.\label{fig: Function value}]{%
        \includegraphics[width=0.32\textwidth]{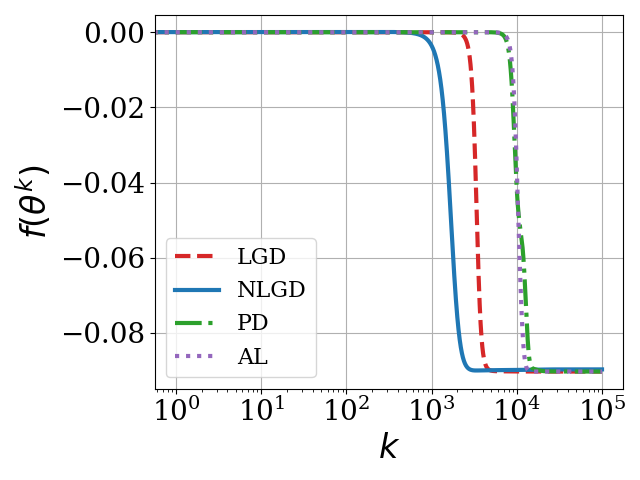}}
    \hfill
    \subfloat[$\|\boldsymbol{\theta}^{k} - \boldsymbol{\theta}^{saddle}\|$.\label{fig: Distance}]{%
        \includegraphics[width=0.32\textwidth]{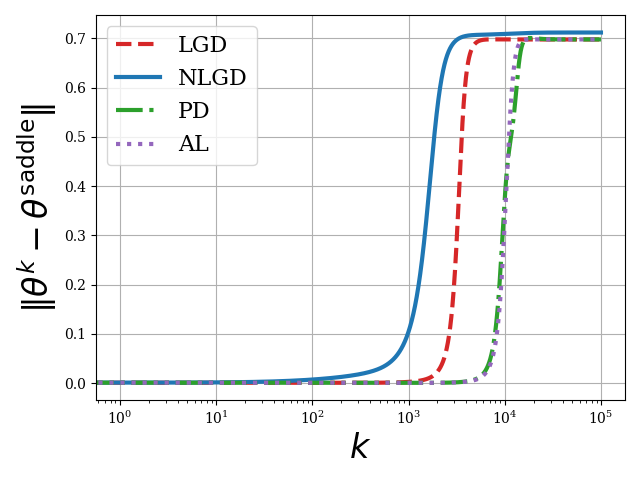}}
    \caption{Second order properties of \textbf{NLGD} and \textbf{LGD} for the smart grid example.}
    \label{fig: Smart grid example}
\end{figure*}

Fig.~\ref{fig: Smart grid example} illustrates the behavior of \textbf{LGD}, \textbf{NLGD}, as well as representative primal--dual (PD) \cite{ding2018primal} and augmented Lagrangian (AL) \cite{ananduta2021distributed} based methods. Fig.~\ref{fig: Function value} shows the evolution of the objective value, while Fig.~\ref{fig: Distance} reports the distance between the iterates and a saddle point.
The results demonstrate that \textbf{NLGD} consistently escapes the saddle point earlier than the \textbf{LGD} baseline and the \textbf{PD}/\textbf{AL} methods, supporting the proposed second-order guarantees.

\begin{figure*}
    \centering
    \subfloat[$f(\boldsymbol{\theta}^{k})$.\label{fig: Function}]{%
        \includegraphics[width=0.32\textwidth]{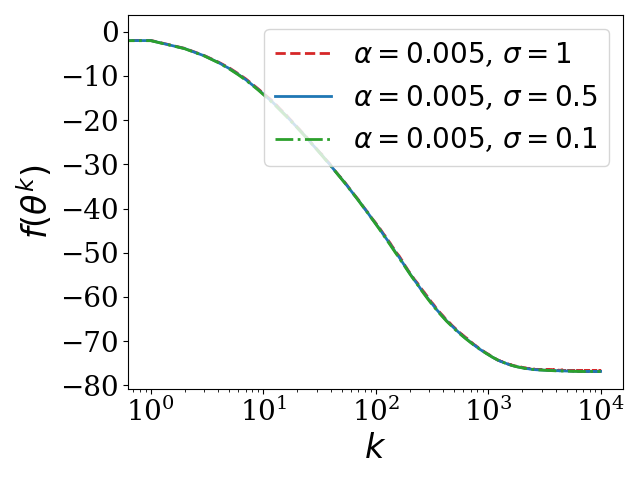}}
    \hfill
    \subfloat[$\|\sqrt{\hat{\mathbf{L}}} \cdot \nabla F(\boldsymbol{\theta}^{k})\|$.\label{fig: Gradient}]{%
        \includegraphics[width=0.32\textwidth]{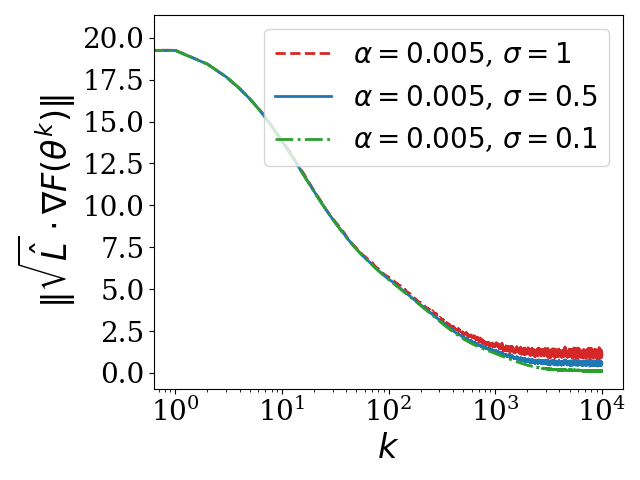}}
    \hfill
    \subfloat[$\min_{\mathbf{d} \in \mathcal{T}} \mathbf{d}^{\top} \nabla^{2} F(\boldsymbol{\theta}^{k}) \mathbf{d} / \|\mathbf{d}\|^{2}$.\label{fig: Hessian}]{%
        \includegraphics[width=0.32\textwidth]{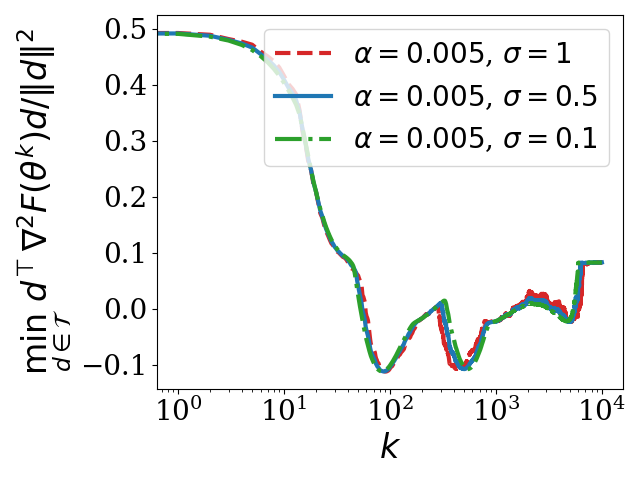}}
    \caption{Second order properties of \textbf{NLGD} for the portfolio example.}
    \label{fig: Portfolio  example}
\end{figure*}

\subsection{Multi-agent Portfolio Optimization}

Another representative example arises in multi-agent portfolio optimization, where multiple decision-making agents (such as institutions or fund managers) allocate capital across a set of financial assets. Each agent aims to balance expected return and risk, while adhering to system-wide constraints such as budget limits or market capacity. A typical objective captures a trade-off between maximizing return and penalizing risk, often using a non-convex regularizer to promote diversification or sparsity in the portfolio. This leads to a distributed non-convex optimization problem of the form
\begin{gather*}
    \min_{\boldsymbol{\theta} \in (\mathbb{R}^{n})^{m}} \quad
    F(\boldsymbol{\theta}) \triangleq \sum_{i=1}^{m} f_{i}(\boldsymbol{\theta}_{i})
    \quad \text{subject to} \quad
    \sum_{i=1}^{m} \boldsymbol{\theta}_{i} = \mathbf{r}.
\end{gather*}

Each agent $i$ controls a portfolio vector $\boldsymbol{\theta}_{i} \in \mathbb{R}^{n}$, and a typical local objective has the form
\begin{gather*}
    f_{i}(\boldsymbol{\theta}_{i})
    = -\mu_{i}^{\top} \boldsymbol{\theta}_{i}
    + \lambda_{i} \boldsymbol{\theta}_{i}^{\top} \Sigma_{i} \boldsymbol{\theta}_{i}
    + \gamma_{i} \log(1 + \|\boldsymbol{\theta}_{i}\|^{2}),
\end{gather*}
where $\mu_{i} \in \mathbb{R}^{n}$ denotes expected returns, $\Sigma_{i} \in \mathbb{R}^{n \times n}$ is a covariance matrix, and $\gamma_i > 0$ weights the non-convex regularizer.

In the portfolio optimization experiment, we simulate a distributed setting with $m = 100$ agents and $n = 5$ assets. All algorithms are initialized using random feasible initializations that satisfy the global budget constraint, without any bias toward saddle points. Communication is governed by the same connected Watts--Strogatz network used in the smart grid example. We implement \textbf{NLGD} with stepsize $\alpha = 0.005$ and vary the noise level $\sigma \in \{0.1, 0.5, 1\}$.

Fig.~\ref{fig: Portfolio example} shows the performance of \textbf{NLGD} under different noise levels. Fig.~\ref{fig: Function} illustrates objective value decrease; Fig.~\ref{fig: Gradient} reports the projected gradient norm
$\|\sqrt{\hat{\mathbf{L}}}\nabla F(\boldsymbol{\theta}^k)\|$; and Fig.~\ref{fig: Hessian} tracks the minimum Rayleigh quotient of the Hessian restricted to the tangent space. The Hessian-based metric becomes negative during the transient phase and subsequently returns to nonnegative values, indicating successful escape from strict saddle points and convergence to second-order stationary points.

\section{Conclusions and Discussion}
\label{sec: Conclusions and Discussion}
This work considers distributed non-convex resource allocation under global constraints and applies \textbf{L}aplacian \textbf{G}radient \textbf{D}escent (\textbf{LGD}) and its newly proposed perturbed variant, \textbf{N}oisy \textbf{LGD} (\textbf{NLGD}). We show that \textbf{LGD} corresponds to gradient descent on an auxiliary function and approaches first-order stationarity. To achieve second-order guarantees, \textbf{NLGD} introduces random perturbations and is shown to reach approximate second-order stationary points with high probability. Numerical experiments on smart grid and portfolio optimization problems validate the theory, demonstrating that \textbf{NLGD} escapes saddle points more effectively and achieves faster convergence than \textbf{LGD}.

The proposed framework focuses on non-convex distributed resource allocation over undirected networks and provides probabilistic guarantees of reaching approximate second-order stationary points. Future work will explore extensions of the proposed framework to more complex constraint sets, time-varying and directed communication networks, and realistic communication impairments such as delays and packet loss. Addressing these settings requires new algorithmic designs and robustness analyses to preserve saddle-point escape and convergence guarantees under asymmetry, asynchrony, and feasibility constraints.

\bibliographystyle{siamplain}
\bibliography{Reference}

\end{document}